\documentclass{article}

\usepackage{amsmath,amssymb,amsthm,vmargin,overpic}
\usepackage{graphicx}
\usepackage{hyperref}
\usepackage{authblk}

\theoremstyle{theorem}

\newtheorem{theorem}{Theorem}

\theoremstyle{definition}

\newtheorem{remark}[theorem]{Remark}

\numberwithin{equation}{section}
\numberwithin{theorem}{section}

\newcommand{\PI}{$\rm{P}_{\rm I}$\,}
\newcommand{\ee}{{\rm e}}
\newcommand{\ii}{{\rm i}}
\newcommand{\dd}{{\rm d}}

\title{On the Riemann--Hilbert approach to asymptotics of tronqu\'ee solutions of Painlev\'e I}
\author[1]{Alfredo Deaño\thanks{alfredo.deanho@uc3m.es}}
\affil[1]{Department of Mathematics, Universidad Carlos III de Madrid, Spain}

\begin{document}

\maketitle

\begin{abstract}
In this paper, we revisit large variable asymptotic expansions of tronqu\'ee solutions of the Painlev\'e I equation, obtained via the Riemann--Hilbert approach and the method of steepest descent. The explicit construction of an extra local parametrix around the recessive stationary point of the phase function, in terms of complementary error functions, makes it possible to give detailed information about exponential-type contributions beyond the standard Poincar\'e expansions for tronqu\'ee and tritronqu\'ee solutions.
\end{abstract}

\section{Introduction}
The Painlev\'e I differential equation (\PI) 
\begin{equation}\label{eq:PI}
y_{xx}=6y^2+x
\end{equation}
is a second order nonlinear ODE which appears as the first one on the list of the six Painlev\'e differential equations, we refer the reader to \cite[Chapter 32]{NIST:DLMF} or \cite{FIKN_2006} for the complete list. Solutions of \PI have important applications in several areas: for example, in random matrix models related to $2$D quantum gravity, they appear when considering certain double scaling asymptotic limits as the size of the matrices tends to infinity, see for example the works of F. David \cite{David_1991,David_1993}, as well as \cite{DFGZJ_1995} and more recently \cite{MSW_2008,MSW_2009} and \cite{BD_PI,DK_2006}. Also, solutions of \PI are used in the description of critical asymptotics of integrable partial differential equations \cite{BT_2013,DGK_2009}, and in connection with the Schr\"{o}dinger equation with a cubic potential \cite{Masoero_2010}.

The differential equation \eqref{eq:PI} admits a $\mathbb{Z}_5$-symmetry that is essential in the analysis: namely, if $y(x)$ is a solution of \eqref{eq:PI}, then so is 
\begin{equation}\label{eq:symZ5}
    v(x)
    =
    \omega^{-2} y\left(\omega^{-1}x\right), \qquad \omega=\ee^{2\pi\ii /5}.
\end{equation}

For this reason, it is natural to consider sectors of angle $2\pi/5$ in the complex plane:
\begin{equation}
    \Omega_k
    =
    \left\{
    z\in\mathbb{C}:
    \frac{\pi}{5}+\frac{2(k-1)}{5}<\arg\, z<\frac{\pi}{5}+\frac{2k}{5}
    \right\}
    ,\qquad k\in\mathbb{Z},
\end{equation}
and because of \eqref{eq:symZ5}, we can restrict the analysis to the sector $\Omega_3$, around the negative real axis.

Generic solutions of \PI are meromorphic functions in the complex plane, with double poles, and furthermore they are transcendental in the sense that no solution of \PI satisfies a first order algebraic differential equation with rational coefficients in the variable $x$, see \cite[Theorem 1.1 and Theorem 13.1]{GLS}. As a consequence, unlike other Painlev\'e equations, \PI does not admit any family of rational, algebraic or classical special function solutions.
 
Because of the transcendental behavior of generic solutions of \PI, asymptotic analysis is particu\-larly relevant, and this task has been carried out extensively over the last decades, with different methods: (i) WKB techniques, following the pioneering work of Kawai and Takei \cite{KT_1996} and applied to generic solutions of \PI, (ii) multiple scales analysis by Joshi and Kruskal \cite{JK_1988}, (iii) hyperasymptotic expansions obtained from the differential equation \cite{OD_2005}, (iv) isomonodromy deformation approach developed by Kapaev \cite{Kapaev_2004} and also Kapaev and Kitaev \cite{KK_1993}. We also mention the recent paper \cite{vSV_2022}, where the authors connect WKB with the isomonodromy approach for generic solutions of \PI. A related problem where asymptotic analysis is very relevant is the location of the poles of the solutions of \PI in the complex plane, see for example \cite{CCH_2015,Costin_1999}.

The asymptotic behavior of generic solutions of \PI involves elliptic functions, as explained in \cite{KK_1993}, see also \cite{Kapaev_2001}, and pole fields are distributed in the whole complex plane. This situation changes for \emph{tronqu\'ee} (truncated) solutions, which are asymptotically free of poles in certain sectors of $\mathbb{C}$. These solutions were studied by Boutroux in \cite{Boutroux_1913}, although the standard characterization currently used is due to Joshi and Kitaev \cite{JK_2001}. Namely, tronqu\'ee solutions are characterised by the asymptotic behaviour 
\begin{equation}\label{eq:asymp_tronquee}
y_0(x)
\sim \frac{x^{1/2}}{\sqrt{6}}\sum_{n=0}^{\infty}\frac{y_{0,n}}{x^{5n/2}}, \qquad |x|\to\infty,
\end{equation}
in a given sector of angle less than $4\pi/5$. This corresponds to two neighboring sectors $\Omega_k$ as constructed before. Boutroux's  \emph{tritronqu\'ee solution} is  asymptotically free of poles as $|x|\to\infty$ in \emph{all but one} sectors $\Omega_k$ of angle $\frac{2\pi}{5}$, and it is characterised by the previous asymptotic behaviour for $|\arg\, x|<\frac{4\pi}{5}$ (and then in any other sector of angle $\frac{8\pi}{5}$ using \eqref{eq:symZ5}).
Tronqu\'ee solutions are particularly important in applications for example in asymptotic analysis of orthogonal polynomials \cite{DK_2006} and in random matrix models \cite{BD_PI}.

The asymptotic expansion \eqref{eq:asymp_tronquee} is usually called perturbative, since it is obtained as a formal power series in a certain small variable, in this case $x^{-5/2}$. The coefficients $y_{0,n}$ can be calculated by inserting \eqref{eq:asymp_tronquee} into the differential equation \eqref{eq:PI} and identifying terms with equal powers of $x$. This leads to the nonlinear recursion
\begin{equation}\label{eq:recak}
y_{0,n+1}
=
\frac{25n^2-1}{8\sqrt{6}}y_{0,n}
-\frac{1}{2}\sum_{m=1}^{k}y_{0,m} y_{0,n+1-m}, \qquad y_{0,0}=1,
\end{equation}
which determines the coefficients $y_{0,n}$ uniquely. Furthermore, it can be shown from this recursion that the coefficients $y_{0,n}$ grow factorially with $n$, which renders the expansion \eqref{eq:asymp_tronquee} divergent for any fixed value of $x$. 

The asymptotic structure of tronqu\'ee solutions of \PI is actually more involved, since \eqref{eq:asymp_tronquee} actually holds for a whole family of such solutions. This is a one-parameter family of solutions of \PI, and this additional parameter appears in exponentially small contributions, that are beyond all orders and therefore asymptotically negligible compared to \eqref{eq:asymp_tronquee}. A more complete ansatz for the asymptotic expansion of tronqu\'ee solutions of \PI can be given in the form of a transseries expansion as $(-x)\to\infty$:
\begin{equation}\label{eq:transy}
y(x)
\sim 
\left(-\frac{x}{6}\right)^{1/2}
\sum_{n=0}^{\infty}\frac{y_{0,n}}{(-x)^{5n/2}}+
(-x)^{1/2}\sum_{k=1}^{\infty} C^k (-x)^{-5k/8}\ee^{-k A (-x)^{5/4}} \sum_{n\geq 0} \frac{y_{k,n}}{(-x)^{5n/4}}.
\end{equation}

Here $A>0$, and the case $k=0$ recovers the standard asymptotic expansion (algebraic in $-x$). Observe that the asymptotic variable is different when $k=0$ and when $k\geq 1$. A much more involved recursion can be worked out for these coefficients, from \cite[Section 1.2]{GIKM_2012} (see also \cite[Section 2]{OD_2005}). In the case of generic solutions of \PI, a full two parameter transseries is needed, involving also logarithmic terms, see for instance \cite{ASV_2012}.

These formal series that appear accompanying exponential/logarithmic terms are often called non-perturbative contributions, since they do not appear directly when using a standard power series ansatz for the solution. Nevertheless, using the powerful tools of Borel transforms, resurgence and alien calculus, see for example \cite{Aniceto2019,Delabaere_2016}, one can see that their structure is intimately related to the original (perturbative) asymptotic expansion, and also to the fact that this expansion is divergent for fixed values of the asymptotic variable. We refer the reader to \cite{Delabaere_2016} for more details on this structure for \PI, see also \cite{GIKM_2012} for the asymptotic behavior of instantons, as well as \cite{Takei_2007,Takei_2008} for results on instanton solutions for  the \PI hierarchy.

The analysis of transseries expansions for different problems in Mathematics and in Physics has been an active topic of research in the last years, see for example \cite{ASV_2012,Baldino2023} and also \cite{MPP2010,Marino2008,Marino2014} and references therein. In string theory, such exponential contributions are related to instantons (with action given by the parameter $A$ introduced before), and they are relevant in large $N$ duality. One particular instance of this construction is given by the topological expansion of the free energy $F_N$ (or of the partition function $Z_N$) in $N\times N$ unitarily invariant ensembles with potential function $V(z)$, which is typically a polynomial. This topological expansion for $F_N$ has the form of a divergent series in the variable $N^{-2}$, and it holds in the so-called one cut case (where the large $N$ limit distribution of eigenvalues is supported on a single interval). This is in a certain sense the most stable configuration, because it corresponds to the dominant stationary point of the potential, but it there are other stationary points it is possible to conceive different limit configurations of eigenvalues supported on several different arcs (multicut case). In that case, the extra stationary points would give exponentially small contributions to the total free energy, we refer the reader to \cite{MSW_2009} for a detailed analysis.

The main goal of this paper is to understand the structure of the transseries expansion for tronqu\'ee solutions of \PI, from the perspective of the Riemann--Hilbert formulation and the steepest descent method associated to it. In several instances of this type of analysis (e.g. asymptotic analysis of orthogonal polynomials or random matrix theory, see for instance \cite{Deift_book,DKMVZ_1999a,DKMVZ_1999b} and many subsequent papers), exponentially small terms are usually estimated and then discarded in the final result. In the case of \PI, we show that, for tronqu\'ee solutions, exponentially small corrections with respect to tritronqu\'ee solutions can be calculated explicitly, using a local parametrix around the recessive stationary point, in terms of the complementary error function. This construction allows us to compute the exponentially small corrections with greater detail than the results given in \cite{Kapaev_2004}. 

Although the steepest descent method does not lead to easy recursion formulas for the coefficients in the asymptotic expansions for Painlev\'e I functions (compared with direct matching in the \PI differential equation, for example), it does provide a rigorous proof of the general structure of such expansions. The explicit addition of exponential terms in this example establishes a potentially fruitful connection between the field of exponential asymptotics and the steepest descent method applied to Riemann--Hilbert problems.

\section{Painlev\'e I Riemann--Hilbert problem}

An essential property of Painlev\'e differential equations, developed by Jimbo and Miwa \cite{JM_1981} and also studied in the seminal work by Flaschka and Newell \cite{FN_1981}, is the fact that they govern isomonodromic deformations of certain linear $2\times 2$ systems of differential equations. More precisely, one considers a function $\Psi(\lambda,x)$ that satisfies simultaneously the two linear differential systems
\begin{equation}\label{eq:systemAU}
    \frac{\dd}{\dd\lambda} \Psi(\lambda,x)
    =
    A(\lambda,x)
    \Psi(\lambda,x),\qquad
    \frac{\dd}{\dd x} \Psi(\lambda,x)
    =
    U(\lambda,x)
    \Psi(\lambda,x),
\end{equation}
with coefficients $A(\lambda,x)$ and $U(\lambda,x)$ that are rational functions of $\lambda$, in the extended complex plane. The compatibility condition $\Psi_{\lambda x}=\Psi_{x\lambda}$ leads to the relation $A_x+AU=U _{\lambda}+UA$, which produces nonlinear relations between the entries of the matrices $A$ and $U$. This  idea has wide applicability in the field of integrable systems, including the theory of orthogonal polynomials and integrable partial differential equations. 

In the present case, this compatibility between two systems of ordinary differential equations corresponds to an isomonodromic deformation of the $\lambda$ equation, and for specific forms of $A(\lambda,x)$ and $U(\lambda,x)$ it produces solutions in the form of Painlev\'e equations in the entries of the matrices $A$ and $U$. For Painlev\'e I, the coefficients of the system are
\begin{equation}
\begin{aligned}
    A(\lambda,x)
    &=
    \begin{pmatrix}
    -z & 2\lambda^2+2y\lambda+x+2y^2 \\ 
     2(\lambda-y)   & z
    \end{pmatrix},\\
    U(\lambda,x)
    &=
    \begin{pmatrix}
    0 & \lambda+2y\\ 
    1 & 0
    \end{pmatrix},
\end{aligned}
\end{equation}
where
\begin{equation}
    z=y_x,\qquad z_x=6y^2+x,
\end{equation}
so $y(x)$ indeed solves Painlev\'e I. The interested reader can find systems for other Painlev\'e equations in \cite[Appendix C]{JM_1981}, \cite[Chapter 32]{NIST:DLMF} and \cite[Chapter 5]{FIKN_2006}, as well as extra examples in \cite{KH_1999}.

As explained in \cite[Chapter 5]{FIKN_2006}, the isomonodromy method consists of defining canonical solutions $\Psi_k(\lambda,x)$ of \eqref{eq:systemAU} in suitable sectors of the complex plane. These sectors are determined by the singularity structure of the matrix $A(\lambda,x)$, and solutions in adjacent sectors are connected via (constant) Stokes matrices. Using further properties of \eqref{eq:systemAU}, one expresses all monodromy data (Stokes multipliers and connection matrices between different singularities) in terms of a minimal set of scalars. This construction can then be used to set up a Riemann--Hilbert (RH) problem for a matrix $Y(\lambda,x)$, closely related to $\Psi(\lambda,x)$, which is the starting point for a large $|x|$ asymptotic analysis.


We consider the Painlev\'e I Riemann--Hilbert problem, as presented by Kapaev in \cite{Kapaev_2004}. Let $\Psi(\lambda,x)$ be a $2\times 2$ matrix valued function such that
\begin{enumerate}
    \item $\Psi(\lambda,x)$ is analytic in $\mathbb{C}\setminus \gamma_{\Psi}$, where \begin{equation}
        \gamma_{\Psi}
        =
        \rho\cup \left(\cup_{j=-2}^{2} \gamma_k\right), \qquad
        \gamma_k=\{\lambda\in\mathbb{C}:\arg\, z=\tfrac{2k\pi}{5}\},\qquad
        \rho=\{\lambda\in\mathbb{C}:\arg\, \lambda=\pi\}
    \end{equation} is the union of contours shown in Figure \ref{fig:gammaPsi} and $k\in\mathbb{Z}$. Furthermore, within each of the sectors, the function $\Psi(\lambda,x)$ admits continuous extension to the boundary.

\begin{figure}
\centerline{\includegraphics[scale=1]{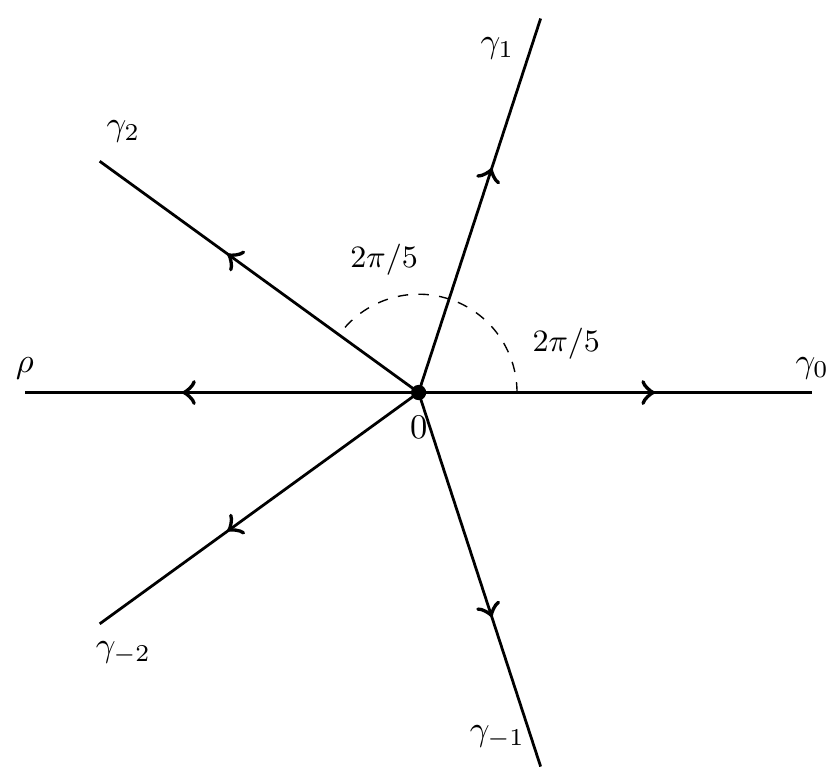}}
\caption{Union of contours $\gamma_{\Psi}$ for the \PI Riemann--Hilbert problem.}
\label{fig:gammaPsi}
\end{figure}

    \item For $\lambda\in\gamma_{\Psi}\setminus\{0\}$, the matrix $\Psi(\lambda,x)$ admits boundary values $\Psi_{\pm}(\lambda,x)$, related by the following jumps:
    \begin{equation}
        \Psi_{+}(\lambda,x)
        =
        \Psi_{-}(\lambda,x)S,
    \end{equation}
    where for $\lambda\in\gamma_k$, the jump matrices are $S=S_k$, with
    \begin{equation}
        S_{2k}=\begin{pmatrix} 1& 0\\ s_{2k} & 1 \end{pmatrix},\qquad
        S_{2k-1}=\begin{pmatrix} 1& s_{2k-1}\\ 0 & 1 \end{pmatrix},
    \end{equation}
    and for $\lambda\in\rho$, we have
    \begin{equation}
        S=\begin{pmatrix} 0 & -\ii\\ -\ii & 0 \end{pmatrix}.
    \end{equation}
    The Stokes multipliers $s_k$, with $k\in\mathbb{Z}$, are constant with respect to $x$, and they verify
    \begin{equation}\label{eq:relations_sk}
        s_{k+5}=s_k, \qquad 1+s_ks_{k+1}=-\ii s_{k+3},\qquad k\in\mathbb{Z}.
    \end{equation}

\item As $\lambda\to\infty$, we have
\begin{equation}\label{eq:Psi_infty}
    \Psi(\lambda,x)
    =
    \lambda^{\frac{\sigma_3}{4}}
    \frac{1}{\sqrt{2}}
    \begin{pmatrix}
    1 & 1\\
    1 & -1
    \end{pmatrix}
    \left(I+\frac{\Psi^{(1)}(x)}{\lambda^{1/2}}+\frac{\Psi^{(2)}(x)}{\lambda}+\mathcal{O}(\lambda^{-3/2})\right)\ee^{\theta_0(\lambda,x)\sigma_3},
\end{equation}
with phase function
\begin{equation}\label{eq:theta0}
\theta_0(\lambda,x)=\frac{4}{5}\lambda^{5/2}+x\lambda^{1/2}.
\end{equation}
The fractional powers in \eqref{eq:Psi_infty} and \eqref{eq:theta0} are taken with their principal branch (real and positive when $\lambda>0$), and with a cut on the contour $\rho$.
\end{enumerate}

The matrices $\Psi^{(1)}(x)$ and $\Psi^{(2)}(x)$ are important in the asymptotic analysis, since they contain functions related to the Painlev\'e I differential equation:
\begin{equation}\label{eq:Psi1Psi2}
    \Psi^{(1)}(x)
    =
    \begin{pmatrix}
    -\mathcal{H}(x) & 0\\
    0 & \mathcal{H}(x) 
    \end{pmatrix},\qquad
    \Psi^{(2)}(x)
    =
    \frac{1}{2}
    \begin{pmatrix}
    \mathcal{H}(x)^2 & y(x)\\
    y(x) & \mathcal{H}(x)^2
    \end{pmatrix},
\end{equation}
where $\mathcal{H}(x)=\frac{z^2}{2}-2y^3-xy$ is the Painlev\'e I Hamiltonian, $z=y_x$ and $z_x=6y^2+x$, so $y(x)$ solves the Painlev\'e I differential equation $y_{xx}=6y^2+x$. The Hamiltonian function  has been a relevant object in the analysis of Painlev\'e equations since the work of Okamoto \cite{Okamoto1980a,Okamoto1980b}, and the related Painlev\'e tau functions appear for example in connection with gap probabilities and expectations of characteristic polynomials in random matrix theory, see for instance \cite{FW_PII}, subsequent papers and references therein.

If we construct the matrix
\begin{equation}
    Y(\lambda,x)
    =
    \frac{1}{\sqrt{2}}
    \begin{pmatrix}
    1 & 1\\ 1 & -1
    \end{pmatrix}
    \lambda^{-\sigma_3/4}
    \Psi(\lambda,x)
    e^{-\theta_0(\lambda,x)\sigma_3},
    \end{equation}
then we can recover the solution to Painlev\'e I as the following limit:
\begin{equation}\label{eq:HfromY}
y(x)
=
2\lim_{\lambda\to\infty} 
\lambda Y_{12}(\lambda,x)
=
2\lim_{\lambda\to\infty} 
\lambda Y_{21}(\lambda,x),
\end{equation}
as well as the Hamiltonian
\begin{equation}\label{eq:yfromY}
\mathcal{H}(x)
=
-\lim_{\lambda\to\infty} 
\lambda^{1/2} Y_{11}(\lambda,x),
\end{equation}
with the principal value of the root and a cut on the negative real axis.

From \eqref{eq:relations_sk}, it follows that fixing two Stokes multipliers will determine all of them. This is consistent with the two degrees of freedom that one expects from solutions of a second order ODE such as \PI. In this formulation, tronqu\'ee and tritronqu\'ee solutions are fixed by certain values of the Stokes multipliers, as shown by Kapaev \cite{Kapaev_2004}.

If we take the Painlev\'e variable in the sector $|\arg (-x)|\leq \frac{2\pi}{5}$, the tronqu\'ee family of solutions is given by $s_0=0$. In this case, there is no jump on the positive real axis in the Riemann--Hilbert problem, and \eqref{eq:relations_sk} implies that
\begin{equation}\label{sk}
s_{2}=s_{-2}=\ii, \qquad s_1+s_{-1}=\ii.
\end{equation}
This leaves one degree of freedom, given by $s_1$ or $s_{-1}$.

The set of solutions of \PI admits other symmetry relations: if $y(x)=f(x,\{s_k\})$ is a solution of \PI with a given set of Stokes multipliers, then
\begin{equation}\label{eq:symPI}
y=\overline{f(\overline{x},\{-\overline{s_{-k}}\})}, \qquad
y=\ee^{\frac{4\pi\ii}{5}} f(\ee^{\frac{2\pi\ii}{5}}x,\{s_{k+2n}\}), \qquad n\in\mathbb{Z}
\end{equation}
are solutions as well.

\section{Main results}
Our first result recovers known asymptotics for tritronqu\'ee Painlev\'e I functions.

\begin{theorem}\label{th:y0H0}
In the reduced case $s_0=s_{-1}=0$, the solution of Painlev\'e I and the correspon\-ding Hamiltonian admit the following asymptotic expansions as $x\to\infty$, for $\varphi=\arg x\in\left[\frac{3\pi}{5},\pi\right]$:
\begin{equation}
\begin{aligned}
    y_0(x)
    &\sim
 \left(-\frac{x}{6}\right)^{1/2}
	\sum_{n=0}^{\infty} \frac{y_{0,n}}{(-x)^{5n/2}},\\
    \mathcal{H}_0(x)
    &\sim
    4\left(-\frac{x}{6}\right)^{3/2}
    \sum_{n=0}^{\infty} \frac{h_{0,n}}{(-x)^{5n/2}}.
\end{aligned}
\end{equation}
Here it is understood that $(-x)^{\alpha}=(\ee^{-\pi\ii} x)^{\alpha}$, and 
the first coefficients are 
\begin{equation}\label{eq:a1b1}
y_{0,0}=1, \qquad y_{0,1}=-\frac{\sqrt{6}}{48}, \qquad h_{0,0}=1, \qquad  h_{0,1}=\frac{\sqrt{6}}{32},
\end{equation}
in agreement with \cite{Kapaev_2004}.
\end{theorem}

The second main result contains the exponentially small corrections of tronqu\'ee solutions with respect to the previous reduced case when $s_{-1}\neq 0$.
\begin{theorem}\label{th:yH}
If $s_0=0$ and $\varphi=\arg x\in\left[\frac{3\pi}{5},\pi\right]$, the tronqu\'ee solutions of Painlev\'e I and the correspon\-ding Hamiltonian admit the following asymptotic expansions as $|x|\to\infty$:
\begin{equation}
\begin{aligned}
    y(x)
    &\sim
    y_0(x)
    +
    (-x)^{1/2}
    \sum_{k=1}^{\infty}
    (-x)^{-5k/8}\ee^{-kA(-x)^{5/4}}
    \sum_{n=0}^{\infty}
    \frac{y_{k,n}}{(-x)^{5n/4}}\\
    \mathcal{H}(x)
    &\sim 
    \mathcal{H}_0(x)
    +
    (-x)^{1/4}
    \sum_{k=1}^{\infty}
    (-x)^{-5k/8}
    \ee^{-kA(-x)^{5/4}}
    \sum_{n=0}^{\infty}
    \frac{h_{k,n}}{(-x)^{5n/4}},
\end{aligned}
\end{equation}
where $y_0(x)$ is the solution of \PI with Stokes multipliers $s_0=s_{-1}=0$ and $s_{\pm 2}=s_{1}=\ii$, and $\mathcal{H}_0(x)$ is the corresponding Hamiltonian function. Here it is understood that $(-x)^{\alpha}=(\ee^{-\pi\ii} x)^{\alpha}$ and 
\begin{equation}\label{eq:A}
    A=\frac{2^{11/4} 3^{1/4}}{5}.%
\end{equation}
The first two coefficients are
\begin{equation}\label{eq:y10h10}
y_{1,0}=\frac{2^{-11/8} 3^{-1/8}s_{-1}}{\sqrt{\pi}},\qquad
h_{1,0}=-\frac{2^{-17/8} 3^{-3/8}s_{-1}}{\sqrt{\pi}},
\end{equation}
in agreement with \cite[Theorem 2.2]{Kapaev_2004}.
\end{theorem}

In the previous theorem, the coefficients $y_{k,n}$ and $h_{k,n}$ can be computed from the information obtained in the steepest descent analysis. This approach is systematic but quite laborious, so a direct approach using \eqref{eq:transy} into the Painlev\'e I ODE (see \cite{GIKM_2012} or \cite{OD_2005}) seems to be more convenient in order to determine these coefficients.

Using the symmetry \eqref{eq:symPI}, we can construct another solution of \PI:
\begin{equation}\label{eq:y1}
y_1(x)=\overline{y_0(\ee^{2\pi\ii}\overline{x},\{-\overline{s_{-k}}\})}.
\end{equation}
This solution corresponds to Stokes multipliers $s_0=s_1=0$ and $s_{\pm 2}=s_{-1}=\ii$, and it gives a result analogous to Theorem \ref{th:yH} in the sector $\varphi=\arg z\in\left[\pi,\frac{7\pi}{5}\right]$, replacing $y_0(x)$ by $y_1(x)$. As a consequence, the difference $y_1(x)-y_0(x)$ is exponentially small in the sector $\left(\frac{3\pi}{5},\frac{7\pi}{5}\right)$, which can be understood as a nonlinear Stokes phenomenon in this context. Furthermore, as shown by Kapaev \cite[(2.72)]{Kapaev_2004}, the solution $y_0(x)$ is actually the tritronqu\'ee one in the larger sector $\left(-\frac{\pi}{5},\frac{7\pi}{5}\right)$.

\section{Steepest descent method}
In this section, we consider $x=|x|\ee^{\ii\varphi}$, with $\varphi=\arg x\in \left[\frac{3\pi}{5},\pi\right]$. We note that, with the exception of the extra local parametrix $\widetilde{P}(\lambda,x)$, the transformations in the steepest descent method essentially appear in \cite[Section 2.1]{Kapaev_2004}, with slight modifications in notation. We redo these steps here to make the paper more self-contained. 

\subsection{Scaling}
We begin with the following scaling:
\begin{equation}\label{eq:PsitoPhi_general}
\Phi(\lambda,x)
=
|x|^{-\frac{\sigma_3}{8}}
\Psi(\lambda |x|^{1/2},x).
\end{equation}

The matrix $\Phi(\lambda,x)$ satisfies the following Riemann--Hilbert problem:
\begin{enumerate}
\item $\Phi(\lambda,x)$ is analytic in $\mathbb{C}\setminus \gamma_{\Phi}$.
\item On $\gamma_{\Phi}\setminus\{0\}$, the matrix $\Phi(\lambda,x)$ has the same jumps as $\Psi(\lambda,x)$.
\item As $\lambda\to\infty$, we have the asymptotic behavior
\begin{equation}\label{eq:Phi_infty_general}
    \Phi(\lambda,x)
    \sim
    \lambda^{\frac{\sigma_3}{4}}
    \frac{1}{\sqrt{2}}
    \begin{pmatrix}
    1 & 1\\
    1 & -1
    \end{pmatrix}
    \left(I+\sum_{k=1}^{\infty}
    \frac{\Psi^{(k)}(x)}{\lambda^{k/2}|x|^{k/4}}\right)\ee^{|x|^{5/4}\theta(\lambda)\sigma_3},\qquad
    \theta(\lambda)=\frac{4}{5}\lambda^{5/2}+\ee^{\ii\varphi}\lambda^{1/2}.
\end{equation}
\end{enumerate}
In the sequel we will write $t=|x|^{5/4}$ for the large parameter in the asymptotic analysis.

\subsection{Auxiliary $g$-function}
An essential step in the asymptotic analysis is the introduction of a suitable $g$ function, in order to normalise the Riemann--Hilbert problem at infinity. 

If we just make the transformation 
\[
S(\lambda,t)=\Phi(\lambda,t)\ee^{-t \theta(\lambda)\sigma_3},
\] 
with phase function $\theta(\lambda)$ given in \eqref{eq:Phi_infty_general}, 
then the jump matrices will contain $\ee^{\pm 2t \theta(\lambda)}$
in the upper (with the $+$ sign) or lower (with the $-$ sign) off-diagonal elements.
We observe that the phase function $\theta(\lambda)$ has stationary points
\begin{equation}
\lambda_{\pm}=\pm\frac{\ee^{\frac{\ii(\varphi-\pi)}{2}}}{2}.
\end{equation}
Given that
$\varphi\in \left[\frac{3\pi}{5},\pi\right]$, we have $\arg\lambda_+\in \left[-\frac{\pi}{5},0\right]$ and 
$\arg\lambda_-\in \left[\frac{4\pi}{5},\pi\right]$, and then 
\begin{equation}
\arg\, \theta(\lambda_{+})
\in\left[\frac{\pi}{2},\pi\right],\qquad
\arg\, \theta(\lambda_{-})
\in\left[-\pi,-\frac{\pi}{2}\right].
\end{equation}
It follows that, unless $\varphi=\frac{3\pi}{5}$, the jump matrices for $S(\lambda,x)$ will blow up at $\lambda=\lambda_-$ as $t\to\infty$, because $\textrm{Re}\,\theta(\lambda_-)<0$, so this transformation cannot be performed directly with the phase function $\theta(\lambda)$. 

For this reason, we introduce an auxiliary $g$ function, that behaves as $\theta(\lambda)$ at infinity and is such that the exponential terms appearing in the new jump matrices do not blow up at the stationary points. One such possibility is the following:
\begin{equation}\label{eq:g}
    g(\lambda)
    =
    \frac{4}{5}\left(\lambda+2\lambda_0\right)^{5/2}- 
    4\lambda_0
    \left(\lambda +2\lambda_0\right)^{3/2}, \qquad 
    \lambda_0=
    \frac{\ee^{\frac{\ii(\varphi-\pi)}{2}}}{\sqrt{6}}.
\end{equation}
The fractional powers in \eqref{eq:g} are taken with their principal branch (real and positive for $\lambda+2\lambda_0>0$), with a branch cut along a ray that joins $-\infty$ and $-2\lambda_0$. 

 As $\lambda\to\infty$, this $g$ function satisfies
\begin{equation}\label{eq:gtheta_infty}
g(\lambda)
=
\frac{4}{5}\lambda^{5/2}-6\lambda_0^2\lambda^{1/2}
-4\lambda_0^3\lambda^{-1/2}+\mathcal{O}(\lambda^{-3/2})
=
\theta(\lambda)
-4\lambda_0^3\lambda^{-1/2}+\mathcal{O}(\lambda^{-3/2}), 
\end{equation}
using \eqref{eq:theta0}, so we can rewrite the asymptotic expansion \eqref{eq:Phi_infty_general} as
\begin{equation}\label{eq:Phi_inftyg}
    \Phi(\lambda,t)
    \sim
    \lambda^{\frac{\sigma_3}{4}}
    \frac{1}{\sqrt{2}}
    \begin{pmatrix}
    1 & 1\\
    1 & -1
    \end{pmatrix}
    \left(I+\sum_{k=1}^{\infty}
    \frac{\Phi^{(k)}(t)}{\lambda^{k/2}t^{k/5}}\right)\ee^{tg(\lambda)\sigma_3}.
\end{equation}

\begin{figure}[ht!]
\centerline{ 
 \includegraphics[scale=1]{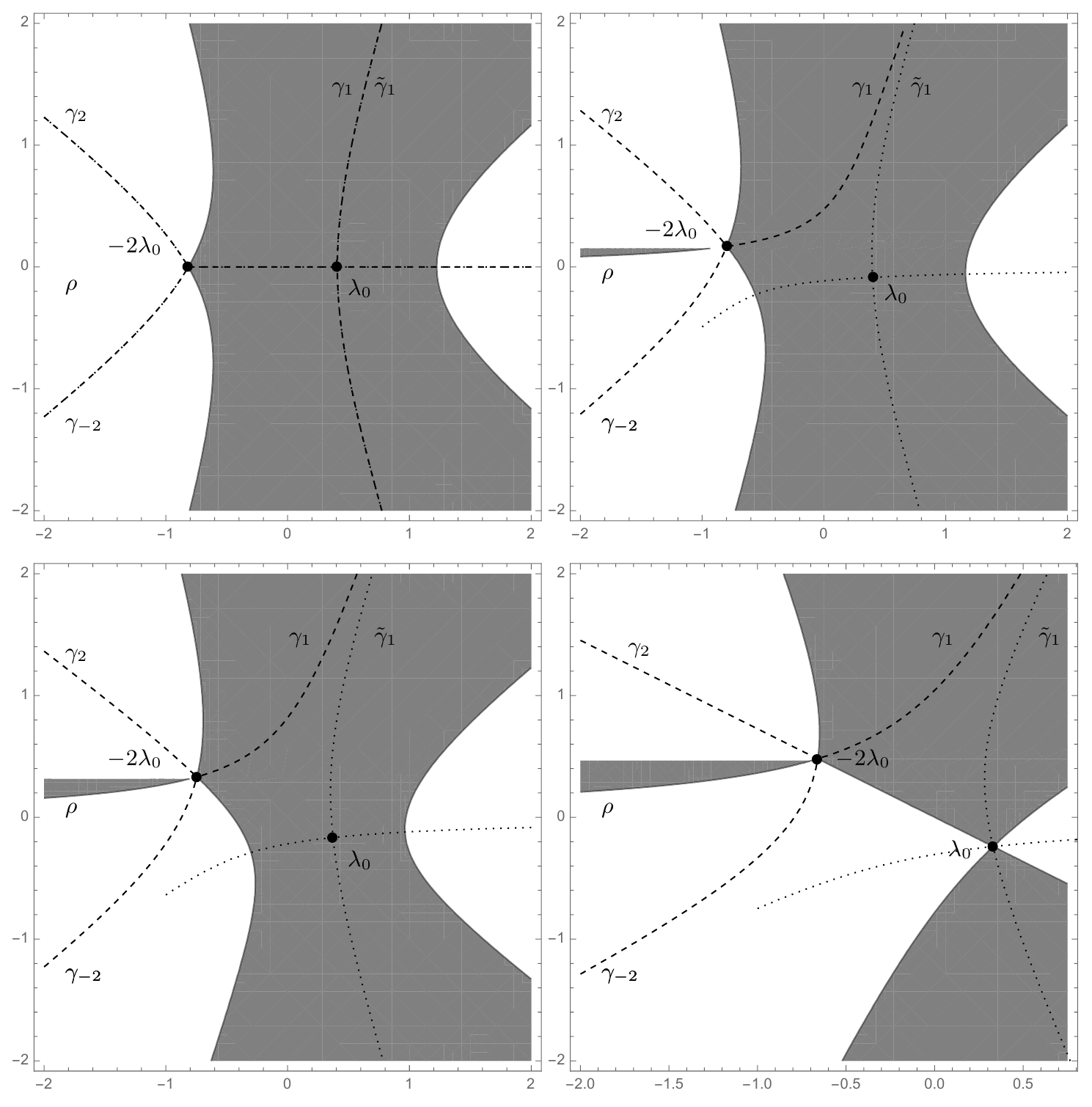}}
\caption{In white, regions where $\textrm{Re}\, g(\lambda)>0$ (white) and where $\textrm{Re}\, g(\lambda)<0$ (gray). From left to right, $\varphi=\pi,\frac{13\pi}{5}$ (first row) and $\varphi=\frac{11\pi}{5},\frac{3\pi}{5}$ (second row). In dashed (respectively dotted) lines, contours of steepest descent through the stationary point $-2\lambda_0$ (respectively $\lambda_0$).}
\label{fig:g0123}
\end{figure}

This modification implies straightforward changes in the coefficients:
\begin{equation}\label{eq:PhiPsi12_general}
\begin{aligned}
    \Phi^{(1)}(t)
    &=
    \Psi^{(1)}(t)
    +4\lambda_0^3 t^{6/5} \sigma_3
    ,\\
    \Phi^{(2)}(t)
    &=
    \Psi^{(2)}(t)
    +
    4\lambda_0^3 t^{6/5}\Psi^{(1)}(t)\sigma_3
    -8\lambda_0^6 t^{12/5} I\\
\end{aligned}
\end{equation}

The function $g(\lambda)$ has two stationary points, which are $-2\lambda_0$ and $\lambda_0$. Note that, for $\varphi\in \left[\frac{3\pi}{5},\pi\right]$, we have 
$\arg \lambda_0\in\left[-\frac{\pi}{5},0\right]$, and $\arg \left(-2\lambda_0\right)\in\left[\frac{4\pi}{5},\pi\right]$. 

We shift the contour 
$\gamma_{\Phi}\mapsto\gamma_{\Phi}-2\lambda_0$, so it is centred at $-2\lambda_0$ instead that at the origin. This implies a redefinition of the matrix $\Phi(\lambda,t)$ in certain sectors of the complex plane, but it does not alter the asymptotics at infinity, so we keep the same notation. We also use the freedom that we have in the original contours of the Riemann--Hilbert problem, because of analyticity, to deform them into curves of steepest descent passing through the stationary points  $-2\lambda_0$ and $\lambda_0$, see Figure \ref{fig:g0123}. The original contour $\rho$ is deformed to the level line $\textrm{Re}\, g(\lambda)=0$, joining $-2\lambda_0$ with $-\infty$.

We note that using the fact that $s_1+s_{-1}=\ii$, the jump on $\gamma_1$ can be factorised as follows:
\begin{equation}
    S_1
    =
    \begin{pmatrix}
    1 & s_1\\
    0 & 1
    \end{pmatrix}
    =
    \begin{pmatrix}
    1 & -s_{-1}\\
    0 & 1
    \end{pmatrix}
    \begin{pmatrix}
    1 & \ii\\
    0 & 1
    \end{pmatrix}
    =
    S_{-1}^{-1}
    \begin{pmatrix}
    1 & \ii\\
    0 & 1
    \end{pmatrix}.
\end{equation}
This will allow us to split the jump on $\gamma_1$ into two jumps. Also, we separate the contour $\gamma_{\Phi}$ into two parts. For this, we need to redefine the matrix $\Phi(\lambda,t)$ in the strip limited by $\gamma_{\pm 1}$ between $-2\lambda_0$ and $\lambda_0$ as follows:
\begin{equation}
\begin{aligned}
   \Phi(\lambda,t)
    &\mapsto
    \Phi(\lambda,t)
    \begin{pmatrix}
    1 & -s_{-1}\\
    0 & 1
    \end{pmatrix}, \qquad \lambda\in\Omega_{\pm},\\
\end{aligned}
\end{equation}

\begin{figure}
\centerline{\includegraphics[scale=1]{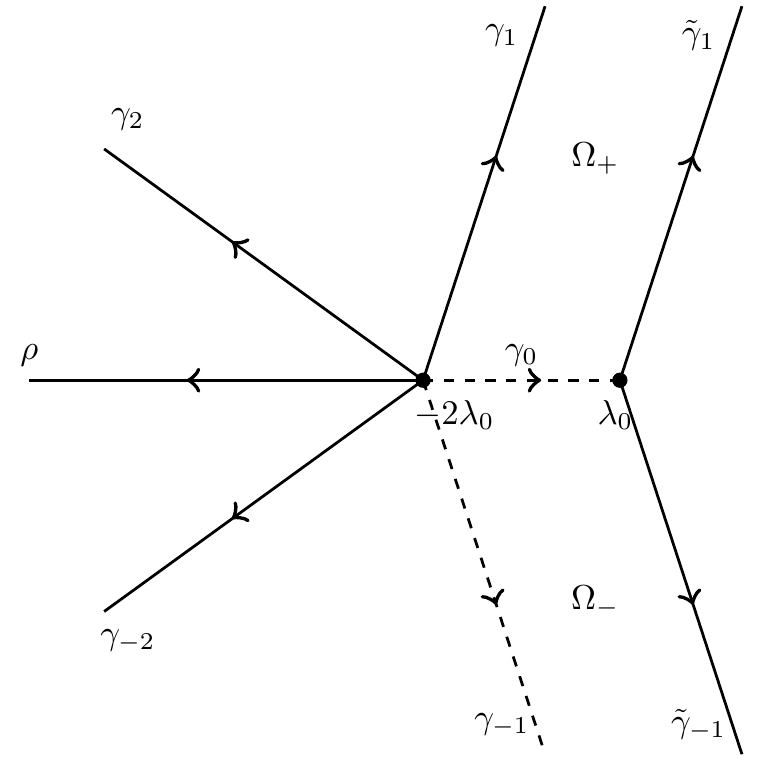}}
\caption{Modified  union of contours $\gamma_{\Phi}$ for the matrix $\Phi(\lambda,t)$. The geometry corresponds to the case $\arg x=\pi$, but it is similar when $\arg x\in\left[\frac{3\pi}{5},\pi\right]$, cf. Figure \ref{fig:g0123}.}
\label{fig:gammaPhit}
\end{figure}

The result is shown in Figure \ref{fig:gammaPhit}. The matrix $\widetilde{\Phi}(\lambda,t)$ satisfies the following Riemann--Hilbert problem:
\begin{enumerate}
\item $\Phi(\lambda,t)$ is analytic in $\mathbb{C}\setminus\gamma_{\widetilde{\Phi}}$, see Figure \ref{fig:gammaPhit}.
\item On $\gamma_{\widetilde{\Phi}}$, the matrix $\widetilde{\Phi}(\lambda,t)$ has the following jumps:
 \begin{equation}
J_{\widetilde{\Phi}}
=
\begin{cases}
\begin{pmatrix}
1 & \ii\\ 0 & 1
\end{pmatrix},&\qquad \lambda\in\gamma_1,\\
\begin{pmatrix}
1 & -s_{-1}\\ 0 & 1
\end{pmatrix},&\qquad \lambda\in\tilde{\gamma}_1,\\
\begin{pmatrix}
1 & s_{-1}\\ 0 & 1
\end{pmatrix},&\qquad \lambda\in\tilde{\gamma}_{-1},\\
\begin{pmatrix}
1 & 0\\ \ii & 1
\end{pmatrix},&\qquad \lambda\in\tilde{\gamma}_{\pm 2},\\
\begin{pmatrix}
0 & -\ii \\ -\ii & 0
\end{pmatrix},&\qquad \lambda\in\rho.
\end{cases}.
\end{equation}
The jumps for the new matrix $\widetilde{\Phi}(\lambda)$ are the same as those of $\Phi(\lambda)$ on $\gamma_{\pm 2}\cup\rho$. 
It is important to notice that there is no jump on the new contour $\gamma_0$, which results from the shift of $\gamma_{\pm 1}$ to $\tilde{\gamma}_{\pm 1}$.
\item As $\lambda\to\infty$, we have the asymptotic behavior
\begin{equation}\label{eq:Phit_inftyg}
    \widetilde{\Phi}(\lambda,t)
    \sim
    (\lambda+2\lambda_0)^{\frac{\sigma_3}{4}}
    \frac{1}{\sqrt{2}}
    \begin{pmatrix}
    1 & 1\\
    1 & -1
    \end{pmatrix}
    \left(I+\sum_{k=1}^{\infty}
    \frac{\widetilde{\Phi}^{(k)}(t)}{\lambda^{k/2}t^{k/5}}\right)\ee^{tg(\lambda)\sigma_3},
\end{equation}
where the coefficients are modified because of the change in the algebraic prefactor:
\begin{equation}\label{eq:PhitPhi12}
\widetilde{\Phi}^{(1)}(t)
=
\Phi^{(1)}(t),\qquad
\widetilde{\Phi}^{(2)}(t)
=
\Phi^{(2)}(t)
-\frac{\lambda_0 t^{2/5}}{2}
\begin{pmatrix} 0 & 1\\ 1 & 0\end{pmatrix}.
\end{equation}

To obtain these modified coefficients $\widetilde{\Phi}^{(1)}(t)$ and $\widetilde{\Phi}^{(2)}(t)$, one possibility is to write the power function in \eqref{eq:Phi_inftyg} as $\lambda^{\frac{\sigma_3}{4}}=(\lambda+2\lambda_0)^{\frac{\sigma_3}{4}}
(\lambda+2\lambda_0)^{-\frac{\sigma_3}{4}}
\lambda^{\frac{\sigma_3}{4}}$, then expand the last two factors as $\lambda\to\infty$ and collect equal powers of $\lambda$.
\end{enumerate}

Next, we make the transformation
\begin{equation}\label{eq:S}
    S(\lambda,t)
    =
    \begin{pmatrix}
    1 & 
    -\widetilde{\Phi}^{(1)}_{11}(t) t^{-1/5}\\
    0 & 1
    \end{pmatrix}
   \widetilde{\Phi}(\lambda,t)
    \ee^{-tg(\lambda)\sigma_3},
\end{equation}
where we recall that $t=|x|^{5/4}$ in terms of the original variable. Then, this new matrix satisfies the following Riemann--Hilbert problem:
\begin{enumerate}
    \item $S(\lambda,t)$ is analytic in $\mathbb{C}\setminus \left(\gamma_{\pm 2}\cup\rho\cup\gamma_1\cup\tilde{\gamma}_{\pm 1}\right)$.
    \item On these contours, the matrix $S$ has the following jumps:
 \begin{equation}\label{eq:jumpsS}
J_{S}
=
\begin{cases}
\begin{pmatrix}
1 & 0 \\ \ii \ee^{-2tg(\lambda)} & 1
\end{pmatrix},&\qquad \lambda\in\gamma_{\pm 2},\\
\begin{pmatrix}
0 & -\ii \\ -\ii & 0
\end{pmatrix},&\qquad \lambda\in\rho,\\
\begin{pmatrix}
1 & \ii\ee^{2tg(\lambda)}\\ 0 & 1
\end{pmatrix},&\qquad \lambda\in\gamma_{1},\\
\begin{pmatrix}
1 & \mp s_{-1}\ee^{2t g(\lambda)}\\ 0 & 1
\end{pmatrix},&\qquad \lambda\in\tilde{\gamma}_{\pm 1},
\end{cases}
\end{equation}
  using the fact that $g_+(\lambda)=-g_-(\lambda)$ for $\lambda\in\rho$. Note that the prefactor introduced in \eqref{eq:S} is multiplying on the left and it is  independent of $\lambda$, so it does not modify any of the jump matrices. 
\item As $\lambda\to\infty$, we have
\begin{equation}\label{eq:S_infty}
    S(\lambda,t)
    =
    \left(I+
    \frac{S^{(1)}(t)}{\lambda}
    +\mathcal{O}(\lambda^{-3/2})\right)
    (\lambda+2\lambda_0)^{\frac{\sigma_3}{4}}
    \frac{1}{\sqrt{2}}
    \begin{pmatrix}
    1 & 1\\
    1 & -1
    \end{pmatrix}.
\end{equation}
This stronger form of the asymptotics (notice that the algebraic factor is now multiplying on the right) holds because of the prefactor that we have introduced when we defined $S(\lambda,t)$. 

\end{enumerate}

We note that
\begin{equation}
    g(-2\lambda_0)=0, \qquad g(\lambda_0)
    =
    	-\frac{24\sqrt{3}}{5}\lambda_0^{5/2}.
\end{equation}
It follows that all jump matrices remain bounded at the two stationary points as $t\to\infty$, and because of the behavior of $\textrm{Re}\, g(\lambda)$ in the complex plane, the jumps on $\gamma_{\pm 2}$ and $\tilde{\gamma}_{\pm 1}$ tend to identity, away from the stationary points.

As a consequence, we can predict that the final contribution from the stationary point $-2\lambda_0$ will be algebraically small in terms of the parameter $t$, whereas the contribution from $\lambda_0$ will be exponentially small as long as $\frac{5}{4}|\varphi-\pi |<\frac{\pi}{2}$, that is $|\varphi-\pi|<\frac{2\pi}{5}$, where $\varphi=\arg x$. In case of equality, the contribution from the stationary point $\lambda=\lambda_0$ will be algebraically small in terms of $t$. This is consistent with the sector for $x$ where we perform the asymptotic analysis, and it establishes the validity of the $g$ function that we have used. 

\subsection{Global parametrix}
We consider the following Riemann--Hilbert problem: we seek $P^{(\infty)}(\lambda)$ such that
\begin{enumerate}
    \item $P^{(\infty)}(\lambda)$ is analytic in $\mathbb{C}\setminus (-\infty,-2\lambda_0]$.
    \item On $(-\infty,-2\lambda_0)$, oriented from left to right, the matrix $P^{(\infty)}(\lambda)$ has the jump
 \begin{equation}\label{eq:jumpPinf}
P^{(\infty)}_+(\lambda)
=
P^{(\infty)}_-(\lambda)
\begin{pmatrix}
0 & \ii \\ \ii & 0
\end{pmatrix}.
\end{equation}
\item As $\lambda\to\infty$, we have
\begin{equation}\label{eq:Pinf_infty}
    P^{(\infty)}(\lambda)
    =
    \left(I+\mathcal{O}(\lambda^{-1/2})\right)
    (\lambda+2\lambda_0)^{\frac{\sigma_3}{4}}
    \frac{1}{\sqrt{2}}
    \begin{pmatrix}
    1 & 1\\
    1 & -1
    \end{pmatrix}.
\end{equation}
\end{enumerate}
The solution of this Riemann--Hilbert problem is given explicitly by
\begin{equation}\label{eq:Pinfty}
    P^{(\infty)}(\lambda)
    =
    (\lambda+2\lambda_0)^{\frac{\sigma_3}{4}}
    \frac{1}{\sqrt{2}}
    \begin{pmatrix}
    1 & 1\\
    1 & -1
    \end{pmatrix},
\end{equation}
where we take the principal branch of the power function, with a branch cut on $(-\infty,-2\lambda_0]$.

Analiticity of this function in $\mathbb{C}\setminus (-\infty,-2\lambda_0]$ and the asymptotic behavior \eqref{eq:Pinf_infty} are straightforward, and the jump condition for $\lambda\in(-\infty,-2\lambda_0)$ can be checked directly as well, using the fact that $(\lambda+2\lambda_0)^{\alpha}_{\pm}
=|\lambda+2\lambda_0|^{\alpha}\ee^{\pm\alpha\pi\ii}$.

\subsection{Local parametrix in a neighborhood of $\lambda=-2\lambda_0$}
Following standard ideas, we take a disc $D(-2\lambda_0,\delta)$ of radius $\delta>0$ fixed around $\lambda=-2\lambda_0$ and we define a matrix $P(\lambda,t)$ that satisfies the following Riemann--Hilbert problem:
\begin{enumerate}
    \item $P(\lambda,t)$ is analytic in $D(-2\lambda_0,\delta)\setminus \left(\gamma_{\pm 2}\cup\rho\cup\gamma_1\right)$.
    \item On these contours, the matrix $P(\lambda,t)$ has the following jumps:
 \begin{equation}
J_{P}
=
\begin{cases}
\begin{pmatrix}
1 & 0 \\ \ii\ee^{-2tg(\lambda)} & 1
\end{pmatrix},&\qquad \lambda\in\gamma_{\pm 2}\cap D(-2\lambda_0,\delta),\\
\begin{pmatrix}
0 & -\ii \\ -\ii & 0
\end{pmatrix},&\qquad \lambda\in\rho\cap D(-2\lambda_0,\delta),\\
\begin{pmatrix}
1 & \ii\ee^{2tg(\lambda)}\\ 0 & 1
\end{pmatrix},&\qquad \lambda\in\gamma_{1}\cap D(-2\lambda_0,\delta).
\end{cases}
\end{equation}
\item Uniformly for $\lambda\in \partial D(-2\lambda_0,\delta)$, we have the matching condition
\begin{equation}\label{eq:P_matching}
    P(\lambda,t)
    =
    P^{(\infty)}(\lambda)(I+\mathcal{O}(t^{-1})), \qquad t\to\infty.
\end{equation}
\end{enumerate}
Now the transformation $\widetilde{P}(\lambda,t)=P(\lambda,t)
\ee^{tg(\lambda)\sigma_3}$
makes constant jumps on the contours inside the disc: 
\begin{enumerate}
    \item $\widetilde{P}(\lambda,t)$ is analytic in $D(-2\lambda_0,\delta)\setminus \left(\gamma_{\pm 2}\cup\rho\cup\gamma_{1}\right)$.
    \item On these contours, the matrix $\widetilde{P}(\lambda,t)$ has the following jumps:
 \begin{equation}
J_{\widetilde{P}}
=
\begin{cases}
\begin{pmatrix}
1 & 0 \\ \ii & 1
\end{pmatrix},&\qquad \lambda\in\gamma_{\pm 2}\cap D(-2\lambda_0,\delta),\\
\begin{pmatrix}
0 & -\ii \\ -\ii & 0
\end{pmatrix},&\qquad \lambda\in\rho\cap D(-2\lambda_0,\delta),\\
\begin{pmatrix}
1 & \ii\\ 0 & 1
\end{pmatrix},&\qquad \lambda\in\gamma_{1}\cap D(-2\lambda_0,\delta).
\end{cases}
\end{equation}
\item Uniformly for $\lambda\in \partial D(-2\lambda_0,\delta)$, we have the matching condition
\begin{equation}\label{eq:Pt_matching}
    \widetilde{P}(\lambda,t)
    =
    P^{(\infty)}(\lambda)(I+\mathcal{O}(t^{-1}))
    \ee^{tg(\lambda)\sigma_3}, \qquad t\to\infty.
\end{equation}
\end{enumerate}

We solve this local problem with Airy functions and a conformal change of variable. Taking into account the local behavior of $g(\lambda)$, namely
\begin{equation}\label{eq:localg}
    g(\lambda)
    =
    -4\lambda_0
    (\lambda+2\lambda_0)^{3/2}
    \left(1+\mathcal{O}(\lambda+2\lambda_0)\right), \qquad \lambda\to 
    -2\lambda_0, 
\end{equation}
with principal branch of the root and a cut along $\rho$, we consider the map
\begin{equation}\label{eq:f}
    f(\lambda)=\left(-\frac{3}{2}g(\lambda)\right)^{2/3}.
\end{equation}    
Using \eqref{eq:localg}, we have the local behavior
\begin{equation}\label{eq:localf}
    f(\lambda)
    =
    (6\lambda_0)^{2/3}(\lambda+2\lambda_0)\left(1+\mathcal{O}(\lambda+2\lambda_0)\right),\qquad \lambda\to -2\lambda_0,
\end{equation}
so $\zeta=f(\lambda)$ is a conformal map from a neighborhood of $\lambda=-2\lambda_0$ onto a neighborhood of $\zeta=0$. Because of the deformation done before to the contours of steepest descent of $g(\lambda)$, then these will be mapped to $\arg\zeta=0$ (in the case of $\gamma_1$), $\arg\zeta=\pm\frac{2\pi}{3}$ (for $\gamma_{\pm 2}$) and $\arg\zeta=\pi$ (in the case of $\rho$).

We seek a local parametrix in the form
\begin{equation}
    \widetilde{P}(\lambda,t)
    =
    E(\lambda,t) A(t^{2/3}f(\lambda))\ee^{-\frac{\pi \ii}{4}\sigma_3},
\end{equation}
where $E(\lambda,t)$ is an analytic prefactor and $A(\zeta)$ will be built out of Airy functions. Namely, we define the following matrix $A(\zeta)$, in the auxiliary $\zeta$ plane, which satisfies the following Riemann--Hilbert problem:
\begin{enumerate}
    \item $A(\zeta)$ is analytic in $\mathbb{C}\setminus \left(\mathbb{R}\cup \ee^{\pm 2\pi \ii/3}[0,\infty)\right)$.
    \item On these rays, oriented towards the origin, $A(\zeta)$ has the following jumps:
    \begin{equation}\label{eq:jumpsA}
    J_{A}(\zeta)
    =
    \begin{cases}	
    \begin{pmatrix}
	1 & 0 \\
	1 & 1
    \end{pmatrix}, \qquad & \arg \zeta=\pm\frac{2\pi}{3},\\
    \begin{pmatrix}
    0 & 1 \\
	-1 & 0
    \end{pmatrix}, \qquad & \zeta\in(-\infty,0)\\
	\begin{pmatrix}
	1 & 1 \\
    0 & 1
    \end{pmatrix}, \qquad & \zeta\in(0,\infty).
	\end{cases}
    \end{equation}
\item As $\zeta\to\infty$, we have the asymptotics 
\begin{equation}
    A(\zeta)
    =
    \frac{1}{2\sqrt{\pi}}\zeta^{-\sigma_3/4}
    \begin{pmatrix}
    1 & \ii\\
    -1 & \ii
    \end{pmatrix}
    \left(I+\mathcal{O}(\zeta^{-3/2})\right)
    \ee^{-\frac{2}{3}\zeta^{3/2}\sigma_3}, \qquad |\arg\,\zeta|<\pi.
\end{equation}
\end{enumerate}
The solution of this Riemann--Hilbert problem is standard, cf. \cite{DKMVZ_1999b}: 
\begin{align}\label{eq: Airy parametrix}
					A(\zeta) = \begin{cases}
					\begin{pmatrix}
					y_0(\zeta) & -y_2(\zeta) \\
					y_0'(\zeta) & -y_2'(\zeta) 
					\end{pmatrix}, \qquad & \arg \zeta \in \left(0, \frac{2\pi}{3}\right), \\
					\begin{pmatrix}
					-y_1(\zeta) & -y_2(\zeta) \\
					-y_1'(\zeta) & -y_2'(\zeta) 
					\end{pmatrix}, \qquad & \arg \zeta \in \left(\frac{2\pi}{3},\pi\right),\\
					\begin{pmatrix}
					-y_2(\zeta) & y_1(\zeta) \\
					-y_2'(\zeta) & y_1'(\zeta) 
					\end{pmatrix}, \qquad & \arg \zeta \in \left(-\pi, -\frac{2\pi}{3}\right),\\
					\begin{pmatrix}
					y_0(\zeta) & y_1(\zeta) \\
					y_0'(\zeta) & y_1'(\zeta) 
					\end{pmatrix}, \qquad & \arg \zeta \in \left(-\frac{2\pi}{3},0\right).
					\end{cases}
				\end{align}
				where 
				\begin{equation}\label{eq: Airy functions}
					y_0(\zeta) := \text{Ai}(\zeta), \qquad y_1(\zeta):=\omega\text{Ai}(\omega \zeta), \qquad y_2(\zeta):=\omega^2 \text{Ai}(\omega^2 \zeta),
				\end{equation}
where $\text{Ai}$ is the Airy function and $\omega:=\exp\left(2\pi\ii/3\right)$. In order to obtain the matching \eqref{eq:Pt_matching} on the boundary of the disc, we need the prefactor
\begin{equation}
    E(\lambda,t)
    =
    P^{(\infty)}(\lambda)
    \ee^{\frac{\pi\ii}{4}\sigma_3}
    \sqrt{\pi}
    \begin{pmatrix}
    1 & -1\\
    -\ii & -\ii
    \end{pmatrix}
    f(\lambda)^{\sigma_3/4}.
\end{equation}
The matrix $E(\lambda,t)$ is analytic in $D(-2\lambda_0,\delta)\setminus(-2\lambda_0-\delta,-2\lambda_0)$. Using the fact that 
\[
f(\lambda)^{\sigma_3/4}_{+}
=\ii\sigma_3
f(\lambda)^{\sigma_3/4}_{-},\qquad 
\lambda\in(-2\lambda_0-\delta,-2\lambda_0),
\]
together with the jump for $P^{(\infty)}(\lambda)$ in \eqref{eq:jumpPinf}, we can show that $E_+(\lambda,t)=E_-(\lambda,t)$ for $\lambda\in(-2\lambda_0-\delta,-2\lambda_0)$. Hence, $E(\lambda,t)$ has analytic continuation across this segment, and the only possible singularity is an isolated one at $\lambda=-2\lambda_0$. By direct calculation we can check that  $E(\lambda,t)=\mathcal{O}(1)$ as $\lambda\to-2\lambda_0$, so $E(\lambda,t)$ is analytic in $D(-2\lambda_0,\delta)$.

Using the classical asymptotic expansions for Airy functions, we can work out the following:
\begin{equation}\label{eq:asympA_full}
A(\zeta)
\sim 
\frac{1}{2\sqrt{\pi}}\zeta^{-\sigma_3/4}
    \begin{pmatrix}
    1 & \ii\\
    -1 & \ii
    \end{pmatrix}
    \left(I+\sum_{k=1}^{\infty}\frac{A_k}{\zeta^{3k/2}}\right)
    \ee^{-\frac{2}{3}\zeta^{3/2}\sigma_3},\qquad |\arg\,\zeta|<\pi,
\end{equation}
where
\begin{equation}\label{eq:Ak}
    A_k
    =
    \frac{(3/2)^k}{2}
    \begin{pmatrix}
    (-1)^k (u_k+v_k) & \ii(u_k-v_k)\\
    (-1)^{k+1} \ii(u_k-v_k) & u_k+v_k
    \end{pmatrix}
\end{equation}
with coefficients $u_0=v_0=1$ and 
\begin{equation}\label{eq:uk}
u_{k}=\frac{(2k+1)(2k+3)(2k+5)\cdots(6k-1)}{216^{k}k!}
=
\frac{\Gamma\left(3k+\frac{1}{2}\right)}{54^k k!\Gamma\left(k+\frac{1}{2}\right)}
,\qquad
v_k=\frac{6k+1}{1-6k}u_k,
\end{equation}
for $k\geq 1$, see \cite[9.7.2]{NIST:DLMF}.

\subsection{Final transformation for $s_0=s_{-1}=0$}
If $s_{-1}=0$, then there is no jump on the contours $\tilde{\gamma}_{\pm 1}$. In this case, we define the following matrix:
\begin{equation}\label{eq:defZ}
Z(\lambda,t)
=
\begin{cases}
    S(\lambda,t)\left(P^{(\infty)}(\lambda)\right)^{-1},&\qquad \lambda\in\mathbb{C}\setminus
    \overline{D(-2\lambda_0,\delta)},\\
    S(\lambda,t)\left(P(\lambda,t)\right)^{-1},&\qquad \lambda\in D(-2\lambda_0,\delta).
\end{cases}
\end{equation}

\begin{figure}
\centerline{\includegraphics[scale=1]{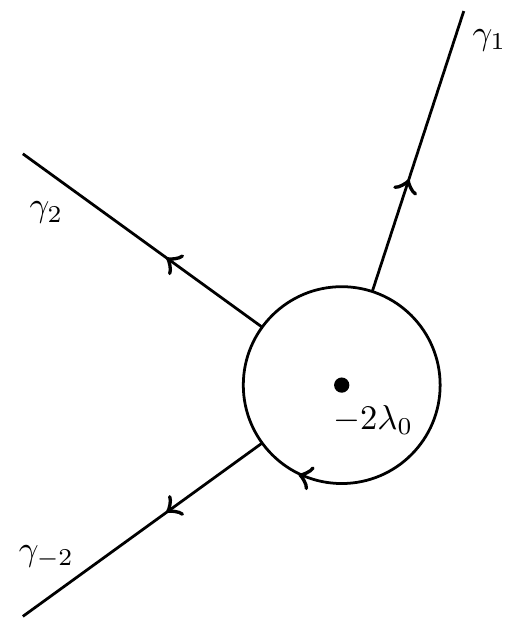}}
\caption{Modified  union of contours $\gamma_Z$ for the final matrix $Z(\lambda,t)$ in the case $s_{-1}=0$.}
\label{figZ}
\end{figure}

Then, $Z(\lambda,t)$ satisfies the following reduced Riemann--Hilbert problem, where we only take into account the jumps on the contours involving the stationary point $-2\lambda_0$:
\begin{enumerate}
    \item $Z(\lambda,t)$ is analytic in $\mathbb{C}\setminus \gamma_Z$, see Figure \ref{figZ}.
    \item For $\lambda\in\gamma_Z$, the matrix $Z(\lambda,t)$ has the following jumps:
    \begin{equation}
        J_{Z}(\lambda,t)
        =
        \begin{cases}
            P^{(\infty)}(\lambda)J_S \left(P^{(\infty)}(\lambda)\right)^{-1},& \qquad
            \lambda\in\left(\gamma_{\pm 2}\cup\gamma_1\right),\\
            P(\lambda,t) \left(P^{(\infty)}(\lambda)\right)^{-1},&
            \qquad
            \lambda\in\partial D(-2\lambda_0,\delta).
        \end{cases}
    \end{equation}
    \item As $\lambda\to\infty$, we have
    \begin{equation}\label{eq:Zinfty}
        Z(\lambda,t)
        =
        I
        +
        \frac{Z^{(1)}(t)}{\lambda}
        +
        \mathcal{O}(\lambda^{-3/2}),
    \end{equation}
with $Z^{(1)}(t)=S^{(1)}(t)$ from before.
\end{enumerate}

On the boundary of the disc $D(-2\lambda_0,\delta)$, we have
\[
J_{Z}(\lambda,t)
=
P(\lambda,t)\left(P^{(\infty)}(\lambda)\right)^{-1}
=
P^{(\infty)}(\lambda)\left(I+\mathcal{O}(t^{-1})\right)
\left(P^{(\infty)}(\lambda)\right)^{-1}
=
I+\mathcal{O}(t^{-1}), \qquad t\to\infty,
\]
because of \eqref{eq:P_matching}.

On $\gamma_{\pm 2}\cup\gamma_1$ outside the disc, the jumps are exponentially close to identity, $\|J_{Z}\|=\mathcal{O}(\ee^{-ct})$, for some constant $c>0$. On the boundary of the disc they are only algebraically close to identity. If we write 
\begin{equation}\label{eq:JZ}
    J_Z(\lambda,t)
    =
    I+\Delta(\lambda,t),
\end{equation}
then we can show that $\Delta(\lambda,t)$ admits an asymptotic expansion in inverse powers of $t$:
\begin{equation}\label{eq:asympDeltat}
    \Delta(\lambda,t)\sim\sum_{k=1}^{\infty}\frac{\Delta_k(\lambda)}{t^k}.
\end{equation}

The coefficients $\Delta_k(\lambda)$ are $0$ for $\lambda\in\gamma_{\pm 2}\cup\gamma_1$ outside of the disc $D(2\lambda_0,\delta)$. On $\partial D(-2\lambda_0,\delta)$, we have
\[
\begin{aligned}
J_{Z}(\lambda,t)
&=
P(\lambda,t)\left(P^{(\infty)}(\lambda)\right)^{-1}\\
&=
E(\lambda,t)A(f(\lambda,t))
\ee^{-\frac{\pi\ii}{4}\sigma_3}
\ee^{-tg(\lambda)\sigma_3}
\left(P^{(\infty)}(\lambda)\right)^{-1}\\
&=
P^{(\infty)}(\lambda)
\ee^{\frac{\pi i}{4}\sigma_3}
\sqrt{\pi}
    \begin{pmatrix}
    1 & -1\\
    -\ii & -\ii
    \end{pmatrix}
    (t^{2/3}f(\lambda))^{\sigma_3/4}
    A(f(\lambda,t))
    \ee^{-\frac{\pi\ii}{4}\sigma_3}
    \ee^{-g(\lambda,t)\sigma_3}
    \left(P^{(\infty)}(\lambda)\right)^{-1}\\
&\sim
P^{(\infty)}(\lambda)
\ee^{\frac{\pi\ii}{4}\sigma_3}
\left(I+\sum_{k=1}^{\infty} \frac{A_k}{t^k f(\lambda)^{3k/2}}\right)
\ee^{-\frac{\pi\ii}{4}\sigma_3}
    \left(P^{(\infty)}(\lambda)\right)^{-1},
\end{aligned}
\]
where the coefficients $A_k$ come from the asymptotic expansion of the Airy functions and are given by \eqref{eq:Ak}. From this calculation, and since $\Delta(\lambda,t)=J_{Z}(\lambda,t)-I$, we deduce that
\begin{equation}\label{eq:DeltakAk}
\Delta_k(\lambda)
=
\frac{1}{f(\lambda)^{3k/2}}
P^{(\infty)}(\lambda)
\ee^{\frac{\pi\ii}{4}\sigma_3}
A_k
\ee^{-\frac{\pi\ii}{4}\sigma_3}
    \left(P^{(\infty)}(\lambda)\right)^{-1},\qquad k\geq 1.
\end{equation}

Direct calculation yields a very simple structure of these coefficients:
\begin{equation}\label{eq:Deltak0}
\begin{aligned}
\Delta_k(\lambda)
&=
\frac{\left(\frac{3}{2}\right)^{k}}
{f(\lambda)^{3k/2}}
\begin{cases}
\begin{pmatrix}
0 & -v_k (\lambda+2\lambda_0)^{1/2}\\
-u_k (\lambda+2\lambda_0)^{-1/2} & 0
\end{pmatrix}
,& \qquad k\,\,\textrm{odd},\\
\begin{pmatrix}
v_k & 0\\
0 & u_k 
\end{pmatrix}
,& \qquad k\,\,\textrm{even}.
\end{cases}
\end{aligned}
\end{equation}

We can check that $\Delta_k(\lambda)$ does not have any jumps on any of the contours inside $D(-2\lambda_0,\delta)$: this follows from  \eqref{eq:Deltak0} and the fact that for $\lambda\in(-\infty,-2\lambda_0)$, we have $f(\lambda)^{3k/2}_+
=
-f(\lambda)^{3k/2}_-$ if $k$ is odd and $f(\lambda)^{3k/2}_+
=
f(\lambda)^{3k/2}_-$ if $k$ is even, and also that $(\lambda+2\lambda_0)^{\alpha}_{\pm}= |\lambda+2\lambda_0|^{\alpha}\ee^{\pm \alpha\pi\ii}$ for $\lambda\in(-2\lambda_0-\delta,-2\lambda_0)$.

Therefore, each coefficient $\Delta_k(\lambda)$ can be continued onto a disc around $-2\lambda_0$ as a meromorphic function, with a pole at $\lambda=-2\lambda_0$. The order of the pole can be computed explicitly, since we have 
\begin{equation}\label{eq:orderpole}
\Delta_{2s-1}(\lambda)
=
\mathcal{O}((\lambda+2\lambda_0)^{-3s+1}),\qquad 
\Delta_{2s}(\lambda)
=
\mathcal{O}((\lambda+2\lambda_0)^{-3s}),\qquad s\geq 1,
\end{equation}
from \eqref{eq:Deltak0} and the local properties of $f(\lambda)$.

The jump matrices for $Z(\lambda,t)$ depend analytically on the parameter $t$, from the previous discussion, and we have the estimates 
\begin{equation}\label{eq:jumpsZ}
\|J_Z(\lambda,t)-I\|
=
\begin{cases}
\mathcal{O}(\ee^{-ct}), &\qquad \lambda\in\gamma_1\cup\gamma_{\pm 2}\,\, \textrm{outside the disc}\,\, D(-2\lambda_0,\delta)\\
\mathcal{O}(t^{-1}), &\qquad \lambda\in\partial D(-2\lambda_0,\delta).
\end{cases}
\end{equation}

The general theory of the steepest descent method yields that the matrix $Z(\lambda,t)$ itself admits a similar asymptotic expansion in inverse powers of $t$:
\begin{equation}\label{eq:asympZt}
Z(\lambda,t)
\sim 
I
+\sum_{j=1}^{\infty}\frac{Z_j(\lambda)}{t^j},\qquad t\to\infty.
\end{equation}

Furthermore, this expansion is uniform in $\lambda\in\mathbb{C}\setminus\partial D(-2\lambda_0,\delta)$, in the following sense: for every $J\geq 1$ there exist constants $C,M$ such that
\begin{equation}
\left\|Z(\lambda,t)-I-\sum_{j=1}^J \frac{Z_j(\lambda)}{t^j}\right\|\leq \frac{C}{|\lambda|t^{J+1}},
\end{equation}
for $|\lambda|>M$. The proof follows from rewriting $Z$ in integral form and using \eqref{eq:JZ} together with \eqref{eq:asympDeltat}. We refer the reader to the seminal paper \cite{DZ_1993}, as well as \cite[Theorem 7.10]{DKMVZ_1999b} for an analogous result in the context of asymptotic analysis of orthogonal polynomials, see also \cite[Theorem 8.1]{FIKN_2006}.

The coefficients $Z_j(\lambda)$ can be computed explicitly, combining \eqref{eq:asympZt} with \eqref{eq:asympDeltat}. More precisely, 
one can write a sequence of additive Riemann--Hilbert problems for the coefficients $Z_j(\lambda)$, on the boundary of the disc $\partial D(-2\lambda_0,\delta)$. Namely, grouping equal powers of $t$ in the jump relation 
$Z_+(\lambda,t)=Z_-(\lambda,t)(I+\Delta(\lambda,t))$, using \eqref{eq:asympDeltat} and \eqref{eq:asympZt}, we obtain
\begin{equation}\label{eq:RHCk}
Z_{j+}(\lambda)
=
Z_{j-}(\lambda)
+
\sum_{m=0}^{j-1} Z_{m-}(\lambda)\Delta_{j-m}(\lambda), \qquad j\geq 1,
\end{equation}
with $Z_0:=I$. 

\subsubsection{The case $j=1$}
The first additive Riemann--Hilbert problem is 
\begin{equation}\label{eq:RHZ1}
Z_{1+}(\lambda)
=
Z_{1-}(\lambda)+\Delta_1(\lambda).
\end{equation}
From \eqref{eq:orderpole}, with $s=1$, we know that $\Delta_1(\lambda)$ is an anti-diagonal matrix and has a pole of order $2$ at $\lambda=-2\lambda_0$, so we write
\begin{equation}
    \Delta_1(\lambda)
    =
    \frac{\Delta^{(-2)}_1}{(\lambda+2\lambda_0)^2}+
    \frac{\Delta^{(-1)}_1}{\lambda+2\lambda_0}
    +
    \mathcal{O}(1)
    ,\qquad \lambda\to-2\lambda_0,
\end{equation}
where the coefficients can be calculated explicitly from \eqref{eq:DeltakAk}:
\begin{equation}
    \Delta^{(-2)}_1
    =
    \frac{1}{288\lambda_0}
    \begin{pmatrix}
    0 & 0\\
    -5 & 0
    \end{pmatrix},\qquad
    \Delta^{(-1)}_1
    =
    \frac{1}{288\lambda_0}
    \begin{pmatrix}
    0 & 7\\
    -\frac{1}{\lambda_0} & 0
    \end{pmatrix}.
\end{equation}

From \eqref{eq:RHZ1}, the matrix $Z_1(\lambda)$ can be computed as a contour integral:
\[
Z_1(\lambda)
=
\frac{1}{2\pi\ii}\oint_{\partial D(-2\lambda_0,\delta)}\frac{\Delta_1(s)}{s-\lambda}\dd s,
\]
with the contour oriented clockwise. Evaluating this integral, we obtain
\begin{equation}\label{eq:Z1}
    Z_1(\lambda)
    =
    \begin{cases}
        \displaystyle
        \frac{\Delta_1^{(-2)}}
        {(\lambda+2\lambda_0)^2}
        +
        \frac{\Delta_1^{(-1)}}{\lambda+2\lambda_0},& \qquad \lambda\in\mathbb{C}\setminus \overline{D(-2\lambda_0,\delta}),\\[3mm]
        \displaystyle
        \frac{\Delta_1^{(-2)}}
        {(\lambda+2\lambda_0)^2}
        +
        \frac{\Delta_1^{(-1)}}{\lambda+2\lambda_0}-\Delta_1
    (\lambda),& \qquad \lambda\in D(-2\lambda_0,\delta).
    \end{cases}
\end{equation}
If we consider the expansion in the outer region, as $\lambda\to\infty$, we obtain
\begin{equation}
\begin{aligned}
    Z_1(\lambda)
    &=
    \frac{1}{288\lambda_0}
    \begin{pmatrix}
    0 & 7\\
    -\frac{1}{\lambda_0} & 0
    \end{pmatrix}
    \frac{1}{\lambda}
    +\mathcal{O}(\lambda^{-2}).
\end{aligned}
\end{equation}

Comparing this with asymptotic expansions \eqref{eq:Zinfty} and \eqref{eq:asympZt}, we can find the first $\mathcal{O}(t^{-1})$ term in the matrix $Z^{(1)}(t)$:
\begin{equation}\label{eq:Zsup1}
Z^{(1)}(t)
=
\frac{1}{288\lambda_0 t}
    \begin{pmatrix}
    0 & 7\\
    -\frac{1}{\lambda_0} & 0
    \end{pmatrix}+\mathcal{O}(t^{-2}), \qquad t\to\infty.
\end{equation}

\subsubsection{The case $j=2$}
The next Riemann--Hilbert problem is given by 
\begin{equation}\label{eq:RHZ2}
Z_{2+}(\lambda)
=
Z_{2-}(\lambda)+Z_{1-}(\lambda)\Delta_1(\lambda)+\Delta_2(\lambda).
\end{equation}
In this case, \eqref{eq:orderpole}, with $s=2$, says that $\Delta_2(\lambda)$ is diagonal and has a pole of order $3$ at $\lambda=-2\lambda_0$, so we write
\begin{equation}
    \Delta_2(\lambda)
    =
    \frac{\Delta^{(-3)}_2}{(\lambda+2\lambda_0)^3}+
    \frac{\Delta^{(-2)}_2}{(\lambda+2\lambda_0)^2}+
    \frac{\Delta^{(-1)}_2}{\lambda+2\lambda_0}
    +
    \mathcal{O}(1)
    ,\qquad \lambda\to-2\lambda_0,
\end{equation}
with
\begin{equation}
\begin{aligned}
    \Delta^{(-3)}_2
    &=
    \frac{1}{2^{11}3^4\lambda_0^2}
    \begin{pmatrix}
    -455 & 0\\
    0 & 385
    \end{pmatrix},\qquad
    \Delta^{(-2)}_2
    =
    \frac{1}{2^{10}3^4\lambda_0^3}
    \begin{pmatrix}
    -91 & 0\\
    0 & 77
    \end{pmatrix},\\
    \Delta^{(-1)}_2
    &=
    \frac{1}{2^{11}3^4 5 \lambda_0^4}
    \begin{pmatrix}
    -91 & 0\\
    0 & 77
    \end{pmatrix}.
\end{aligned}
\end{equation}
With a calculation similar to the previous case, we have 
\[
Z_2(\lambda)
=
\frac{1}{2\pi\ii}\oint_{\partial D(-2\lambda_0,\delta)}\frac{Z_{1-}(s)\Delta_1(s)+\Delta_2(s)}{s-\lambda}\dd s,
\]
with the contour oriented clockwise. Evaluating this integral, we obtain
obtain
\begin{equation}\label{eq:Z2}
    Z_2(\lambda)
    =
    \begin{cases}
        \displaystyle
        \frac{Z^{(-3)}_2}{(\lambda+2\lambda_0)^3}+
    \frac{Z^{(-2)}_2}{(\lambda+2\lambda_0)^2}+
    \frac{Z^{(-1)}_2}{\lambda+2\lambda_0},& \quad \lambda\in\mathbb{C}\setminus \overline{D(-2\lambda_0,\delta}),\\[3mm]
        \displaystyle
        \frac{Z^{(-3)}_2}{(\lambda+2\lambda_0)^3}+
    \frac{Z^{(-2)}_2}{(\lambda+2\lambda_0)^2}+
    \frac{Z^{(-1)}_2}{\lambda+2\lambda_0}-Z_{1-}(\lambda)\Delta_1(\lambda)-\Delta_2(\lambda),& \quad \lambda\in D(-2\lambda_0,\delta).
    \end{cases}
\end{equation}
Here
\begin{equation}
\begin{aligned}
Z_2^{(-3)}&=\frac{35}{2^{11}3^4\lambda_0^2}\begin{pmatrix}-13 &0\\ 0 & 11\end{pmatrix},\quad
Z_2^{(-2)}=\frac{7}{2^{10}3^4\lambda_0^3} \begin{pmatrix}-12 &0\\ 0 & 11\end{pmatrix},\quad
Z_2^{(-1)}=\frac{49}{2^{11}3^4\lambda_0^4} \begin{pmatrix}-1 &0\\ 0 & 1\end{pmatrix}.
\end{aligned}
\end{equation}
If we expand as $\lambda\to\infty$, we obtain
\begin{equation}\label{eq:Z2_asymp}
Z_2(\lambda)
=
\frac{49}{2^{11}3^4\lambda_0^4} \begin{pmatrix}-1 &0\\ 0 & 1\end{pmatrix}\frac{1}{\lambda}+\mathcal{O}(\lambda^{-2}),\qquad \lambda\to\infty.
\end{equation}
Comparing this with the asymptotic expansions \eqref{eq:Zinfty} and \eqref{eq:asympZt} again, we can improve \eqref{eq:Zsup1} to
\begin{equation}\label{eq:Zsup2}
Z^{(1)}(t)
=
\frac{1}{288\lambda_0 t}
    \begin{pmatrix}
    0 & 7\\
    -\frac{1}{\lambda_0} & 0
    \end{pmatrix}
    +
    \frac{49}{2^{11}3^4\lambda_0^4 t^2} \begin{pmatrix}-1 &0\\ 0 & 1\end{pmatrix}+\mathcal{O}(t^{-3}),\qquad t\to\infty.
\end{equation}

\begin{remark}\label{rem:Z}
In principle, we could compute higher order coefficients by residue calculation, but this is not a very efficient method. We make, however, two important observations: i) we can prove that $Z_{m}(\lambda)$ is anti-diagonal if $m$ is odd and it is diagonal if $m$ is even, which follows from this same property for the matrices $\Delta_{m}(\lambda)$ and \eqref{eq:RHCk}; ii) this approach proves rigorously the existence of a full asymptotic expansion in inverse powers of $t=|x|^{5/4}$ for the final matrix $Z(\lambda,t)$, which translates into a similar expansion for the matrix $Z^{(1)}(t)$ and consequently for the Painlev\'e I functions.
\end{remark}

\begin{remark}\label{rem:expsmall_2l0}
On the contours $\gamma_1\cup\gamma_{\pm 2}$ outside the disc $D(-2\lambda_0,\delta)$, the jump matrices are exponentially close to identity. In standard steepest descent analysis, this estimation is enough to discard them in comparison with the jump on the boundary of the disc, which is only algebraically close to identity. However, when $s_{-1}\neq 0$, we want to add exponentially small contributions from the other stationary point $\lambda=\lambda_0$, and for this reason we have to be particularly careful. Namely, it is not enough to show that $c>0$ in the estimate $\|J_Z(\lambda,t)-I\|=\mathcal{O}(\ee^{-ct})$ for $\lambda\in\gamma_1\cup\gamma_{\pm 2}$ outside the disc $D(-2\lambda_0,\delta)$, but we also need to ensure that $c$ can be made arbitrarily large, in such a way that these contributions are still smaller than the exponential contributions that we will add later on. 

\begin{figure}[h!]
\centerline{\includegraphics[scale=1]{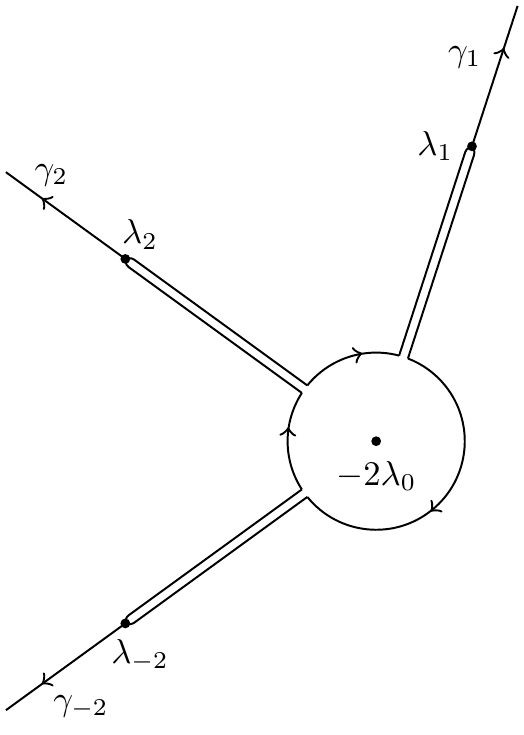}}
\caption{Deformed contour for the final matrix $Z(\lambda,t)$ in the case $s_{-1}=0$.}
\label{fig_blob}
\end{figure}

Using \eqref{eq:defZ} and \eqref{eq:jumpsS}, the jump for $Z(\lambda,t)$ on $\gamma_{\pm 2}$ outside of the disc is
\[
\begin{aligned}
J_Z(\lambda,t)
=
P^{(\infty}(\lambda)J_S(\lambda,t)\left(P^{(\infty}(\lambda)\right)^{-1}
&=
P^{(\infty}(\lambda)\begin{pmatrix} 1 &0\\ \ii\ee^{-2tg(\lambda)} & 1\end{pmatrix}\left(P^{(\infty)}(\lambda)\right)^{-1}\\
&=
I+\ii\ee^{-2t g(\lambda)}P^{(\infty)}(\lambda)\begin{pmatrix} 0 &0\\ 1 & 0\end{pmatrix}\left(P^{(\infty)}(\lambda)\right)^{-1}
\end{aligned}
\]
and on $\gamma_{1}$ outside of the disc is
\[
\begin{aligned}
J_Z(\lambda,t)
=
P^{(\infty}(\lambda)J_S(\lambda,t)\left(P^{(\infty}(\lambda)\right)^{-1}
&=
P^{(\infty}(\lambda)\begin{pmatrix} 1 & \ii\ee^{2tg(\lambda)}\\ 0 & 1\end{pmatrix}\left(P^{(\infty)}(\lambda)\right)^{-1}\\
&=
I+\ii\ee^{2tg(\lambda)}P^{(\infty)}(\lambda)\begin{pmatrix} 0 & 1\\ 0 & 0\end{pmatrix}\left(P^{(\infty)}(\lambda)\right)^{-1}.
\end{aligned}
\]

In order to make these jumps arbitrarily small in $t$, we deform the boundary of the disc $D(-2\lambda_0,\delta)$ following the paths of steepest descent/ascent of $g(\lambda)$, as shown in Figure \ref{fig_blob}. By analiticity, this deformation does not alter the contribution from the stationary point. We denote by $\lambda_{\pm 2}$ and $\lambda_1$ the points at the tips of these extended legs of the contour, and we define
\[
M:=\min \{ \textrm{Re}\,g(\lambda_2),\textrm{Re}\,g(\lambda_{-2}),-\textrm{Re}\,g(\lambda_1)\}.
\] 

Because the contours $\gamma_{\pm 2}$ and $\gamma_1$ follow the paths of steepest descent of $g(\lambda)$, if we move along them far away enough from $-2\lambda_0$, given any $c_0>0$, we can make $M$ as large (and positive) as we wish. In particular, if we take $M=c_0/2$, we obtain
\[
\begin{aligned}
\|J_Z(\lambda,t)-I\|
<
C \ee^{-2t M}=C \ee^{-tc_0},
\end{aligned}
\]
for some $C>0$, since the global parametrix $P^{(\infty)}(\lambda)$ is bounded in any norm for $|\lambda+2\lambda_0|>\delta>0$. With an argument similar to \cite[Theorem 8.1]{FIKN_2006}, this allows us to deduce that the contributions to the large $t$ asymptotic expansion of $Z(\lambda,t)$ coming from the contours away from the disc are arbitrarily small, even compared with any exponential term like the ones that we will be adding later on.
\end{remark}

\subsection{Local parametrix in a neighborhood of $\lambda=\lambda_0$}
In order to study tronqu\'ee solutions, we need to take into account the jumps on the contours $\tilde{\gamma}_{\pm 1}$, see Figure \ref{fig:gammaPhit}, which are not equal to identity if $s_{-1}\neq 0$. In the work of Kapaev \cite[Theorems 2.2 and 2.3]{Kapaev_2004}, an estimate of the size of the jumps on these two contours is enough to establish the nonlinear Stokes phenomenon and the exponentially small corrections with respect to the tritronqu\'ee solution. However, it turns out that we can work out an explicit local parametrix in a neighborhood of $\lambda_0$, and that allows us to provide more detailed information.

We note that, as $\lambda\to\lambda_0$, we have
\begin{equation}\label{eq:glocal_l0}
    g(\lambda)
    =
    g_0
    +
    \sqrt{3\lambda_0}(\lambda-\lambda_0)^2
    +
    \mathcal{O}((\lambda-\lambda_0)^3),
\qquad
    g_0
    =
    -\frac{24\sqrt{3}}{5}\lambda_0^{5/2}.
\end{equation}
Note that $\arg\lambda_0\in\left[-\frac{\pi}{5},0\right]$, so $\arg\lambda_0^{5/2}\in\left[-\frac{\pi}{2},0\right]$, and then $\textrm{Re}\, g_0<0$ in the sector that we are considering.

We take a disc $D(\lambda_0,\delta)$ and we 
define a matrix $Q$ that satisfies the following Riemann--Hilbert problem:
\begin{enumerate}
    \item $Q(\lambda,t)$ is analytic in $D(\lambda_0,\delta)\setminus 
    \tilde{\gamma}_{\pm 1}$.
    \item On these contours, the matrix $Q$ has the following jumps:
 \begin{equation}
J_{Q}
=
\begin{pmatrix}
1 & \pm s_{-1} \ee^{2tg(\lambda)}\\ 0 & 1
\end{pmatrix},\qquad \lambda\in\tilde{\gamma}_{\pm 1}\cap D(\lambda_0,\delta).
\end{equation}
\item Uniformly for $\lambda\in \partial D(\lambda_0,\delta)$, we have the matching condition
\begin{equation}\label{eq:Q_matching}
    Q(\lambda,t)
    =
    P^{(\infty)}(\lambda)
    \left(I+\mathcal{O}(t^{-1/2})\right), \qquad t\to\infty.
\end{equation}
\end{enumerate}

As before, the transformation $\widetilde{Q}(\lambda,t)=Q(\lambda,t)\ee^{t g(\lambda)\sigma_3}$ makes constant jumps on the two contours $\tilde{\gamma}_{\pm 1}\cap D(\lambda_0,\delta)$.

We look for this local parametrix in the following form:
\begin{equation}
    \widetilde{Q}(\lambda,t)
    =
    \widetilde{E}(\lambda)
    B\left(t^{1/2}(g(\lambda)-g_0)^{1/2}\right)
\end{equation}
We note that $\xi=(g(\lambda)-g_0)^{1/2}$ is a conformal map from a neighborhood of $\lambda=\lambda_0$ onto a neighborhood of $\xi=0$, because of the local behavior \eqref{eq:glocal_l0}. The root is taken with its principal value (real and positive when $g(\lambda)-g_0>0$).

Let $B(\xi)$ satisfy the following model Riemann--Hilbert problem:
\begin{enumerate}
    \item $B(\xi)$ is analytic in $\mathbb{C}\setminus \ii\mathbb{R}$.
    \item For $\xi\in \ii\mathbb{R}$, oriented upwards, the matrix $B(\xi)$ has the jump
    \begin{equation}
        B_+(\xi)=B_-(\xi)
        \begin{pmatrix} 
        1 & -s_{-1} \\ 0 & 1 
        \end{pmatrix}.
    \end{equation}
    \item As $\xi\to\infty$, we have
    \begin{equation}
        B(\xi)
        =
        \left(I+\mathcal{O}(\xi^{-1})\right)\ee^{\xi^2\sigma_3}.
    \end{equation}
\end{enumerate}

This Riemann-Hilbert problem can be solved explicitly in terms of the complementary error function (as it was  done in local analysis of orthogonal polynomials, see for example  \cite{BCD_2021,BDY_2017}). If we define
\begin{equation}
    b(\xi)=\frac{s_{-1} \ee^{\xi^2}}{2}
    \begin{cases}
    -\textrm{erfc}(-\sqrt{2}\xi), & \textrm{Re}\,\xi<0,\\
    \textrm{erfc}(\sqrt{2}\xi), & \textrm{Re}\,\xi>0,
    \end{cases}
\end{equation}
then using the property $\textrm{erfc}(z)+\textrm{erfc}(-z)=2$, see \cite[7.4.2]{NIST:DLMF}, we have 
\[
b_+(\xi)-b_-(\xi)
=
-\frac{s_{-1} \ee^{\xi^2}}{2}
\left(\textrm{erfc}(-\sqrt{2}\xi)+
\textrm{erfc}(\sqrt{2}\xi)\right)
=
-s_{-1} \ee^{\xi^2}.
\]
As a consequence, the matrix
\[
B(\xi)
=
\begin{pmatrix}
\ee^{\xi^2} & b(\xi)\\
0 & \ee^{-\xi^2}
\end{pmatrix}
=
\begin{pmatrix}
1 & \ee^{\xi^2}b(\xi)\\
0 & 1
\end{pmatrix}
\ee^{\xi^2\sigma_3}
\]
satisfies the required jump on the imaginary axis. Furthermore, using the asymptotic expansion for the complementary error function 
\cite[7.12.1]{NIST:DLMF}, we have
\begin{equation}\label{eq:asymp_B}
    B(\xi)
    \sim 
    \left(I+\sum_{k=0}^{\infty} 
    \frac{B_k}{\xi^{2k+1}}\right)\ee^{\xi^2\sigma_3},
    \qquad
    B_k
    =\begin{pmatrix}
    0 & b_k\\
    0 & 0
    \end{pmatrix},
    \qquad b_k=\frac{(-1)^k}{2^{k+\frac{3}{2}}\sqrt{\pi}} \left(\frac{1}{2}\right)_k
    s_{-1}.
\end{equation}
for $\xi\to\infty$ with $|\arg \xi|<\frac{3\pi}{4}$. It then follows that 
\[
\begin{aligned}
\ee^{tg_0\sigma_3}
B\left(t^{1/2}(g(\lambda)-g_0)^{1/2}\right)
&=
\ee^{tg_0\sigma_3}
\left(I+\sum_{k=0}^{\infty}\frac{B_k}{t^{k+\frac{1}{2}}(g(\lambda)-g_0)^{k+\frac{1}{2}}}
\right)
    \ee^{t(g(\lambda)-g_0)\sigma_3}\\
&=
\left(I
+
\ee^{2tg_0}
\sum_{k=0}^{\infty}
\frac{B_k}{t^{k+\frac{1}{2}}(g(\lambda)-g_0)^{k+\frac{1}{2}}}
\right)
    \ee^{t g(\lambda)\sigma_3}.
\end{aligned}
\]
In order to obtain the matching \eqref{eq:Q_matching} on the boundary of the disc, we take the following prefactor:
\begin{equation}\label{eq:Et}
\widetilde{E}(\lambda,t)
=
P^{(\infty)}(\lambda)
\ee^{tg_0\sigma_3},
\end{equation}
which is analytic in $D(\lambda_0,\delta)$. 

As a consequence,
\[
\begin{aligned}
Q(\lambda,t)
&=
\widetilde{E}(\lambda,t)
B\left(t^{1/2}(g(\lambda)-g_0)^{1/2}\right)
\ee^{-t g(\lambda)\sigma_3}
=
P^{(\infty)}(\lambda)
\left(
I
+
\ee^{2tg_0}\sum_{k=0}^{\infty}
\frac{B_k}{t^{k+\frac{1}{2}}(g(\lambda)-g_0)^{k+\frac{1}{2}}}
\right)
\end{aligned}
\]
as $t\to\infty$, uniformly for $\lambda\in\partial D(\lambda_0,\delta)$.

\subsection{Final transformation for $s_{-1}\neq 0$}
The final matrix $R(\lambda,t)$ in the steepest descent method is defined as follows:
\begin{equation}\label{eq:defR}
    R(\lambda,t)
    =
   \begin{cases}
    S(\lambda,t)\left(P^{(\infty)}(\lambda)\right)^{-1},&\qquad \lambda\in\mathbb{C}\setminus\left(\overline{D(-2\lambda_0,\delta)}\cup \overline{D(\lambda_0,\delta)}\right),\\
    S(\lambda,t)\left(P(\lambda,t)\right)^{-1},&\qquad \lambda\in D(-2\lambda_0,\delta),\\
    S(\lambda,t)\left(Q(\lambda,t)\right)^{-1},& \qquad \lambda\in D(\lambda_0,\delta).
    \end{cases}
\end{equation}

\begin{figure}
\centerline{\includegraphics[scale=1]{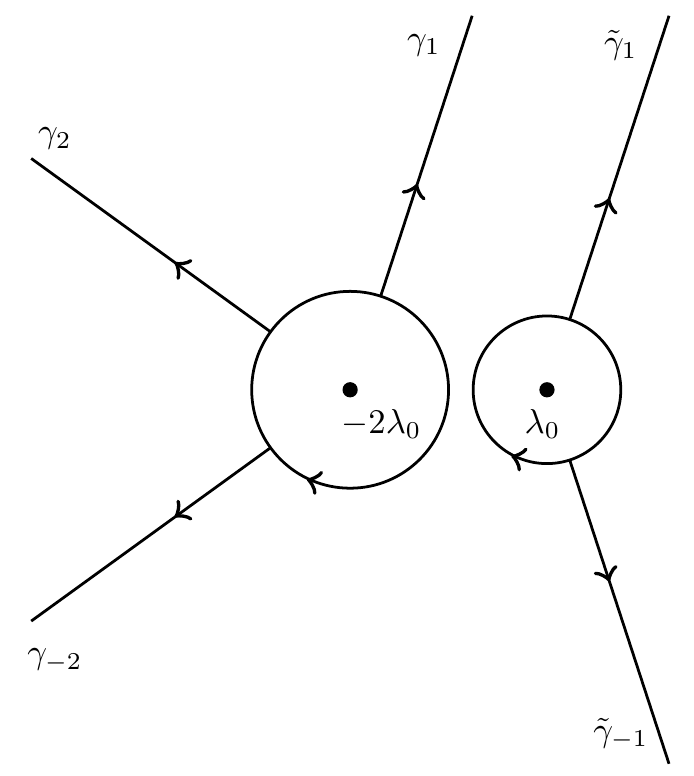}}
\caption{Modified  union of contours $\gamma_R$ for the final matrix $R(\lambda,x)$.}
\label{figR}
\end{figure}

It satisfies the following Riemann--Hilbert problem:
\begin{enumerate}
    \item $R(\lambda,t)$ is analytic in $\mathbb{C}\setminus \gamma_R$, see Figure \ref{figR}.
    \item For $\lambda\in\gamma_R$, the matrix $R(\lambda)$ has the following jumps:
    \begin{equation}
        J_R(\lambda,x)
        =
        \begin{cases}
            P^{(\infty)}(\lambda)J_S \left(P^{(\infty)}(\lambda)\right)^{-1},& \qquad
            \lambda\in\left(\gamma_{\pm 2}\cup\gamma_1\cup\tilde{\gamma}_{\pm 1}\right),\\
            P(\lambda,t) \left(P^{(\infty)}(\lambda)\right)^{-1},&
            \qquad
            \lambda\in\partial D(-2\lambda_0,\delta),\\
            Q(\lambda,t)
            \left(P^{(\infty)}(\lambda)\right)^{-1},&
            \qquad
            \lambda\in\partial D(\lambda_0,\delta).
        \end{cases}
    \end{equation}
    \item As $\lambda\to\infty$, we have
    \begin{equation}\label{eq:Rinfty}
        R(\lambda,t)
        =
        I+
        \frac{R^{(1)}(t)}{\lambda}
        +\mathcal{O}(\lambda^{-2}).
    \end{equation}
\end{enumerate}

Let us consider the jumps of $R(\lambda,t)$ on the different contours:
\begin{itemize}
    \item On $\gamma_{\pm 2}$, outside of the disc $D(-2\lambda_0,\delta)$, we have
\[
J_{R}(\lambda,t)
=
P^{(\infty)}(\lambda)
\begin{pmatrix}
1 & 0\\ \ii\ee^{-2tg(\lambda)} & 1
\end{pmatrix}
\left(P^{(\infty)}(\lambda)\right)^{-1}
=
I
+
\ii\ee^{-2tg(\lambda)}
P^{(\infty)}(\lambda)
\begin{pmatrix}
0 & 0\\ 1 & 0
\end{pmatrix}
\left(P^{(\infty)}(\lambda)\right)^{-1}.
\]

\item On $\gamma_{1}$ outside of the disc $D(-2\lambda_0,\delta)$, we have
\[
J_{R}(\lambda,t)
=
P^{(\infty)}(\lambda)
\begin{pmatrix}
1 & \ii\ee^{2tg(\lambda)}\\ 0 & 1
\end{pmatrix}
\left(P^{(\infty)}(\lambda)\right)^{-1}
=
I
+
\ii\ee^{2tg(\lambda)}
P^{(\infty)}(\lambda)
\begin{pmatrix}
0 & 1\\ 0 & 0
\end{pmatrix}
\left(P^{(\infty)}(\lambda)\right)^{-1}.
\]

\item On $\partial D(-2\lambda_0,\delta)$, we have
\[
J_{R}(\lambda,t)
=
P(\lambda,t)\left(P^{(\infty)}(\lambda)\right)^{-1}
=
P^{(\infty)}(\lambda)\left(I+\mathcal{O}(t^{-1})\right)
\left(P^{(\infty)}(\lambda)\right)^{-1}
=
I+\mathcal{O}(t^{-1}),
\] 
as $t\to\infty$,
because of the matching \eqref{eq:P_matching}.
\end{itemize}

In order to study the jumps on $\partial D(\lambda_0,\delta)\cup \tilde{\gamma}_{\pm 1}$, we write the final matrix $R(\lambda,t)$ as product of two: 
\begin{equation}\label{eq:RCchi}
    R(\lambda,t)
    =
    \chi(\lambda,t)Z(\lambda,t),
\end{equation}
where $Z(\lambda,t)$ is the solution of the reduced Riemann--Hilbert problem 
($s_0=s_{-1}=0$) presented before. The correction matrix $\chi(\lambda,t)=R(\lambda,t)Z(\lambda,t)^{-1}$ then satisfies the following Riemann--Hilbert problem:

\begin{enumerate}
    \item $\chi(\lambda,t)$ is analytic in $\mathbb{C}\setminus \left(\partial D(\lambda_0,\delta)\cup \tilde{\gamma}_{\pm 1}\right)$.
    \item For $\lambda\in\partial D(\lambda_0,\delta)\cup \tilde{\gamma}_{\pm 1}$, the matrix $\chi(\lambda,t)$ has the following jumps:
    \begin{equation}\label{eq:jumps_chi}
        J_{\chi}(\lambda,t)
       =
        \begin{cases}
           Z(\lambda,t)Q(\lambda,t)        \left(P^{(\infty)}(\lambda)\right)^{-1}
            Z(\lambda,t)^{-1},&
            \qquad
            \lambda\in\partial D(\lambda_0,\delta),\\
            Z(\lambda,t)P^{(\infty)}(\lambda)J_S
            \left(P^{(\infty)}(\lambda)\right)^{-1}
            Z(\lambda,t)^{-1},&
            \qquad
            \lambda\in\tilde{\gamma}_{\pm 1}.
        \end{cases}
    \end{equation}
    \item As $\lambda\to\infty$, we have
    \begin{equation}\label{eq:asympchi}
        \chi(\lambda,t)=I+\frac{\chi^{(1)}(t)}{\lambda}+\mathcal{O}(\lambda^{-2}).
    \end{equation}
\end{enumerate}

The jump $J_{\chi}(\lambda,t)$ coincides with $J_R(\lambda,t)$ on $\partial D(\lambda_0,\delta)\cup \tilde{\gamma}_{\pm 1}$, by construction. 
Similarly as before, we write the jump matrix as 
\begin{equation}\label{eq:Jchi}
J_{\chi}(\lambda,t)=I+\widehat{\Delta}(\lambda,t),
\end{equation}
and we show that $\widehat{\Delta}(\lambda,t)$ admits an asymptotic expansion in inverse powers of $t$ as $t\to\infty$, but also incorporating an exponentially small factor:
\begin{equation}\label{eq:asympDeltahat}
    \widehat{\Delta}(\lambda,t)
    \sim
    t^{1/2}\ee^{2t g_0}
    \sum_{j=1}^{\infty}\frac{\widehat{\Delta}_j(\lambda)}{t^{j}}, \qquad  
    g_0
    =
    -\frac{24\sqrt{3}}{5}\lambda_0^{5/2}.
\end{equation}
On $\tilde{\gamma}_{\pm 1}$ outside of the disc $D(\lambda,\delta)$, the jump matrix is 
\[
\begin{aligned}
J_{\chi}(\lambda,t)
&=
Z(\lambda,t)P^{(\infty)}(\lambda)
\begin{pmatrix}
1 & \mp s_{-1}\ee^{2tg(\lambda)}\\ 0 & 1
\end{pmatrix}
\left(P^{(\infty)}(\lambda)\right)^{-1}Z(\lambda,t)^{-1}\\
&=
I
\mp s_{-1}\ee^{2tg(\lambda)}
Z(\lambda,t)P^{(\infty)}(\lambda)
\begin{pmatrix}
0 & 1\\ 0 & 0
\end{pmatrix}
\left(P^{(\infty)}(\lambda)\right)^{-1}Z(\lambda,t)^{-1}.
\end{aligned}
\]

Along the path of steepest descent through the stationary point $\lambda_0$, we have
\[
\textrm{Re}\,g(\lambda)
\leq 
\textrm{Re}\,g_0
=
-\frac{2^{7/4}3^{1/4}}{5}
\cos\left(\tfrac{5}{4}(\varphi-\pi)\right)
\]
and therefore, the jump matrix on $\tilde{\gamma}_{\pm 1}$ outside of the disc $D(\lambda,\delta)$ is exponentially close to identity as $t\to\infty$, $\|J_{\chi}(\lambda,t)-I\|=\mathcal{O}(\ee^{-\hat{c} t})$, with $\hat{c}>2\,|\textrm{Re}\,g_0|$. Then, we have $\widehat{\Delta}_k(\lambda)=0$ on $\tilde{\gamma}_{\pm 1}$ in formula \eqref{eq:Jchi}, and furthermore, with a deformation like the one indicated in Remark \ref{rem:expsmall_2l0}, we can make the constant $\hat{c}$ arbitrarily large, and so negligible with respect to any contribution from the boundary of the disc $ \partial D(\lambda_0,\delta)$.

 On $\partial D(\lambda_0,\delta)$, we have
\[
\begin{aligned}
J_{\chi}(\lambda)
&=
Z(\lambda,t)Q(\lambda)
\left(P^{(\infty)}(\lambda)\right)^{-1}
Z(\lambda,t)^{-1}\\
&\sim
Z(\lambda,t)
P^{(\infty)}(\lambda)
\left(I
+
\ee^{2tg_0}
\sum_{j=0}^{\infty}
\frac{B_j}{t^{j+\frac{1}{2}}(g(\lambda)-g_0)^{j+\frac{1}{2}}}
    \right)
\left(P^{(\infty)}(\lambda)\right)^{-1}
Z(\lambda,t)^{-1}\\
&\sim 
I
+
t^{1/2}\ee^{2tg_0}
\sum_{j=1}^{\infty} 
\frac{\widehat{\Delta}_j(\lambda)}{t^{j}},
\end{aligned}
\]
because of \eqref{eq:Q_matching}. The matrices $B_j$ come from the asymptotic expansion of the local parametrix and they are given in \eqref{eq:asymp_B}. This yields the coefficients $\widehat{\Delta}_j(\lambda,t)$ as a combination of the expansion of the local parametrix around $\lambda=\lambda_0$ and the terms $\Delta_j(\lambda)$ that appear in the asymptotic expansion for $Z(\lambda,t)$, see \eqref{eq:asympZt}.

The first ones are
\begin{equation}\label{eq:Dt1Dt2}
\begin{aligned}
    \widehat{\Delta}_1(\lambda)
    &=
    \frac{1}{(g(\lambda)-g_0)^{1/2}}
    P^{(\infty)}(\lambda)
    B_0
    \left(P^{(\infty)}(\lambda)\right)^{-1},\\
    \widehat{\Delta}_2(\lambda)
    &=
    \frac{1}{(g(\lambda)-g_0)^{1/2}}
    \left[
    Z_1(\lambda),P^{(\infty)}(\lambda)B_0\left(P^{(\infty)}(\lambda)\right)^{-1}
    \right]
    +
    \frac{1}{(g(\lambda)-g_0)^{3/2}}
	P^{(\infty)}(\lambda) B_1
    \left(P^{(\infty)}(\lambda)\right)^{-1},
\end{aligned}
\end{equation}
where $[A,B]=AB-BA$ is the standard matrix commutator. 

Given the local behavior of $g(\lambda)-g_0$ as $\lambda\to \lambda_0$, it is straightforward to check that $\widehat{\Delta}_k(\lambda)$ does not have any jump inside $D(\lambda_0,\delta)$, and that it can be extended to a meromorphic function with a pole of order $1$ at $\lambda=\lambda_0$. Note that the order of the pole does not change with $k$, because in principle there is always at least one non-zero term multiplying the factor $(g(\lambda)-g_0)^{-1/2}$.

As before, the jump matrices for $\chi(\lambda,t)$ depend analytically on the parameter $t$, and then the matrix $\chi(\lambda,t)$ itself admits a similar asymptotic expansion in inverse powers of $t$. However, because of the exponential term $\ee^{2tg_0}$, we have to make a more complicated ansatz, with two sequences of coefficients:

\begin{equation}\label{eq:asympchit}
\chi(\lambda,t)
\sim
I+t\sum_{k=1}^{\infty}\frac{\ee^{2kt g_0}}{t^{\frac{k}{2}}} \sum_{j=1}^{\infty}\frac{\chi_{k,j}(\lambda)}{t^{j}}.
\end{equation}
In particular, note that the expansion for $\chi(\lambda,t)$ needs to have an infinite sequence of exponential terms, even if the jump has only one. Furthermore, this large $t$ asymptotic expansion is uniform in $\lambda\in\mathbb{C}\setminus \partial D(\lambda_0,\delta)$, and $\chi_{k,j}(\lambda)=\mathcal{O}(\lambda^{-1})$ as $\lambda\to\infty$.

If we combine the equation
\begin{equation}\label{eq:chijump}
\chi_+(\lambda)=\chi_-(\lambda)(I+\widehat{\Delta}(\lambda,t)),
\end{equation}
with \eqref{eq:asympchit} and \eqref{eq:asympDeltahat}, we obtain
\begin{equation}\label{eq:addRHchi}
\begin{aligned}
 I+t\sum_{k=1}^{\infty}\frac{\ee^{2kt g_0}}{t^{\frac{k}{2}}} \sum_{j=1}^{\infty}\frac{\chi_{k,j+}(\lambda)}{t^{j}}
 =
 \left(
 I+t\sum_{k=1}^{\infty}\frac{\ee^{2kt g_0}}{t^{\frac{k}{2}}} \sum_{j=1}^{\infty}\frac{\chi_{k,j-}(\lambda)}{t^{j}}
 \right)
 \left(I+t^{1/2}\ee^{2t g_0}
     \sum_{j=1}^{\infty}\frac{\widehat{\Delta}_j(\lambda)}{t^{j}}\right).
\end{aligned}
\end{equation}

At each exponential level (i.e. for each value of $k$), we can construct a sequence of additive Riemann--Hilbert problems for the terms in the large $t$ asymptotic expansion for $\chi(\lambda)$.

\subsubsection{The case $k=1$}
For $k=1$, if we collect terms multiplying the exponential term $\ee^{2tg_0}$, we obtain
\begin{equation}
\chi_{1,j+}(\lambda)
=
\chi_{1,j-}(\lambda)+\widehat{\Delta}_j(\lambda),\qquad j\geq 1.
\end{equation}

The first term, corresponding to $j=1$, can be calculated directly using \eqref{eq:Dt1Dt2} and $b_0=2^{-3/2}s_{-1}\pi^{-1/2}$ from \eqref{eq:asymp_B}. This gives
\begin{equation}
\widehat{\Delta}_1(\lambda)
=
\frac{\widehat{\Delta}_1^{(-1)}}{\lambda-\lambda_0}+\mathcal{O}(1),\qquad 
    \widehat{\Delta}_1^{(-1)}
    =
    \frac{2^{-5/2}s_{-1}}{\sqrt{\pi}(3\lambda_0)^{1/4}}
    \begin{pmatrix}
    1 & -\sqrt{3\lambda_0}\\
    1/\sqrt{3\lambda_0} & -1\\
    \end{pmatrix}
\end{equation}
as $\lambda\to\lambda_0$. From this, we can compute the first correction matrix explicitly:
\begin{equation}
\chi_{1,1}(\lambda)
=
\begin{cases}
        \displaystyle
        \frac{\widehat{\Delta}_1^{(-1)}}
        {\lambda-\lambda_0},& \qquad \lambda\in\mathbb{C}\setminus \overline{D(\lambda_0,\delta}),\\[3mm]
        \displaystyle
        \frac{\widehat{\Delta}_1^{(-1)}}
        {\lambda-\lambda_0}-\widehat{\Delta}_1(\lambda),& \qquad \lambda\in D(\lambda_0,\delta).
\end{cases}
\end{equation}
If we expand as $\lambda\to\infty$, we have
\begin{equation}
\chi_{1,1}(x)
=
\frac{\widehat{\Delta}_1^{(-1)}}{\lambda}+\mathcal{O}(\lambda^{-2}),
\end{equation}
and comparing this with \eqref{eq:asympchi} and the expansion \eqref{eq:asympchit}, we obtain
\begin{equation}\label{eq:chi111}
\chi^{(1)}_{1,1}(x)
=
\ee^{2tg_0}t^{-1/2}\left(\widehat{\Delta}_1^{(-1)}+\mathcal{O}(t^{-1})\right).
\end{equation}


\subsubsection{The case $k\geq 2$}
Equating terms that multiply $\ee^{2ktg_0}$ in  \eqref{eq:addRHchi}, for $k\geq 2$, we obtain
\begin{equation}
\chi_{k,j+}(\lambda)
=
\chi_{k,j-}(\lambda)
+
\sum_{\ell=1}^{j}
\chi_{k-1,\ell-}(\lambda)
\widehat{\Delta}_{j+1-\ell}(\lambda),\qquad j\geq 1.
\end{equation}

These Riemann--Hilbert problems can be solved sequentially, using the matrices $\widehat{\Delta}_j(\lambda)$ and the necessary $\chi_{k-1,j}(\lambda)$ terms from the previous steps. The leading term for each exponential can be calculated using residues:
\[
\chi_{k,1}(\lambda)
=
\frac{1}{2\pi\ii}\oint_{\partial D(\lambda_0,\delta)}\frac{\chi_{k-1,1-}(u)\widehat{\Delta}_1(u)}{u-\lambda}\dd u,
\]
with the boundary of the disc oriented clockwise.

As before, the method of calculation of coefficients can be complicated, particularly for higher values of $k$ and $j$. The procedure is however systematic.

\section{Proof of Theorems \ref{th:y0H0} and \ref{th:yH}}
Now we consider the original matrix $Y(\lambda,x)$, and from \eqref{eq:Psi1Psi2}, we have
\begin{equation}\label{eq:H0y0fromY}
\mathcal{H}(x)
=
-\lim_{\lambda\to\infty} \lambda^{1/2}Y_{11}(\lambda,x)\qquad
y_0(x)
=
2\lim_{\lambda\to\infty} \lambda Y_{12}(\lambda,x).
\end{equation}
We undo the transformations of the steepest descent method, using that $t=|x|^{5/4}$:
\begin{equation}
\begin{aligned}
Y(\lambda,x)
&=
\frac{1}{\sqrt{2}}\begin{pmatrix} 1 & 1\\ 1 & -1 \end{pmatrix}
\lambda^{-\frac{\sigma_3}{4}}\Psi(\lambda,x)\ee^{-\theta_0(\lambda,x)\sigma_3}\\
&=
\frac{1}{\sqrt{2}}\begin{pmatrix} 1 & 1\\ 1 & -1 \end{pmatrix}
\lambda^{-\frac{\sigma_3}{4}}|x|^{\frac{\sigma_3}{8}}\Phi(\lambda |x|^{-1/2},x)
\ee^{-|x|^{5/4}\theta(\lambda |x|^{-1/2})\sigma_3}\\
&=
\frac{1}{\sqrt{2}}\begin{pmatrix} 1 & 1\\ 1 & -1 \end{pmatrix}
\lambda^{-\frac{\sigma_3}{4}}|x|^{\frac{\sigma_3}{8}}
\begin{pmatrix}
1 & \widetilde{\Phi}_{11}^{(1)}(x)|x|^{-1/4}\\ 0 & 1
\end{pmatrix}
S(\lambda |x|^{-1/2},x)\\
&\times 
\ee^{|x|^{5/4}\left(g(\lambda|x|^{-1/2})-\theta(\lambda |x|^{-1/2})\right)\sigma_3}\\
&=
\frac{1}{\sqrt{2}}\begin{pmatrix} 1 & 1\\ 1 & -1 \end{pmatrix}
\lambda^{-\frac{\sigma_3}{4}}|x|^{\frac{\sigma_3}{8}}
\begin{pmatrix}
1 & \widetilde{\Phi}_{11}^{(1)}(x)|x|^{-1/4}\\ 0 & 1
\end{pmatrix}
R(\lambda |x|^{-1/2},x)P^{(\infty)}(\lambda|x|^{-1/2})\\
&\times
\ee^{|x|^{5/4}\left(g(\lambda|x|^{-1/2})-\theta(\lambda |x|^{-1/2})\right)\sigma_3}.
\end{aligned}
\end{equation}

Using the asymptotic expansions \eqref{eq:Rinfty} and \eqref{eq:Pinfty}, we obtain
\begin{equation}
\widetilde{\Phi}_{11}^{(1)}(x)-R_{21}^{(1)}(x)|x|^{1/4}=0,
\end{equation}
and using this relation, we get
\begin{equation}\label{eq:HYR}
\mathcal{H}(x)
=
-\lim_{\lambda\to\infty} \lambda^{1/2}Y_{11}(\lambda,x)
=
4\lambda_0^3|x|^{3/2}-R_{21}^{(1)}(x)|x|^{1/4}
\end{equation}
and
\begin{equation}\label{eq:yYR}
y(x)
=
2\lim_{\lambda\to\infty} \lambda Y_{12}(\lambda,x)
=
\lambda_0 |x|^{1/2}+\left(\left(R_{21}^{(1)}(x)\right)^2+R_{11}^{(1)}-R_{22}^{(1)}\right) |x|^{1/2}. 
\end{equation}

If $s_0=s_{-1}=0$, then $R=Z$, and then
\begin{equation}\label{eq:H0y0Z}
\begin{aligned}
\mathcal{H}_0(x)
&=
4\lambda_0^3|x|^{3/2}-Z_{21}^{(1)}(x)|x|^{1/4}\\
y_0(x)
&=
\lambda_0 |x|^{1/2}+\left(\left(Z_{21}^{(1)}(x)\right)^2+Z_{11}^{(1)}(x)-Z_{22}^{(1)}(x)\right) |x|^{1/2}. 
\end{aligned}
\end{equation}
Since $Z$ admits an asymptotic expansion in inverse powers of $t=|x|^{5/4}$, and this expansion is uniform in $\lambda$, the same is true for $Z^{(1)}$, and therefore for the Painlev\'e functions. Furthermore, because of the structure of the matrices $Z_m$ (see Remark \ref{rem:Z}) and the specific combination of entries that appears in \eqref{eq:H0y0Z}, the asymptotic expansion for $y_0(x)$ is in fact an expansion in inverse powers of $t^2=|x|^{5/2}$. The first coefficient can be computed from \eqref{eq:Zsup2}, namely
\[
\begin{aligned}
\left(\left(Z_{21}^{(1)}(x)\right)^2+Z_{11}^{(1)}(x)-Z_{22}^{(1)}(x)\right)|x|^{1/2}
&=
\lambda_0 |x|^{1/2}\left(-\frac{1}{1728\lambda_0^5 |x|^{5/2}}+\mathcal{O}(|x|^{-5})\right)\\
&=
-\frac{1}{48 x^2}+\mathcal{O}(x^{-9/2}),
\end{aligned}
\]
which gives $y_{0,1}=-\frac{\sqrt{6}}{48}$ and coincides with well known results and also with \eqref{eq:recak}.

In the case $s_{-1}\neq 0$, we have $R(\lambda,x)=\chi(\lambda,x)Z(\lambda,x)$, and from \eqref{eq:Rinfty} we obtain $R^{(1)}(x)=\chi^{(1)}(x)+Z^{(1)}(x)$. Therefore, from \eqref{eq:yYR} and \eqref{eq:HYR}, 
\begin{equation}\label{eq:HyZ}
\begin{aligned}
\mathcal{H}(x)
&=
4\lambda_0^3|x|^{3/2}-(\chi_{21}^{(1)}(x)+Z_{21}^{(1)}(x))|x|^{1/4},\\
&=
\mathcal{H}_0(x)-\chi_{21}^{(1)}(x)|x|^{1/4}\\
y(x)
&=
\lambda_0 |x|^{1/2}+\left(\left(\chi_{21}^{(1)}(x)+Z_{21}^{(1)}(x)\right)^2+\chi_{11}^{(1)}(x)+Z_{11}^{(1)}(x)-\chi_{22}^{(1)}(x)-Z_{22}^{(1)}(x)\right) |x|^{1/2}\\
&=
y_0(x)+\left(2\chi_{21}^{(1)}(x)Z_{21}^{(1)}(x)+\left(\chi_{21}^{(1)}(x)\right)^2
+\chi_{11}^{(1)}(x)-\chi_{22}^{(1)}(x)\right)|x|^{1/2}. 
\end{aligned}
\end{equation}

At this point, the calculation of the coefficients $y_{k,j}$ and $h_{k,j}$ for the Painlev\'e functions in Theorem \ref{th:yH} relies only on the terms $Z^{(1)}(x)$ and $\chi^{(1)}(x)$ that appear in the steepest descent method. This calculation becomes more and more complicated with increasing $k$ (smaller exponential terms) and $j$ (smaller algebraic terms within each exponential level), but it can be tackled systematically. More importantly, \eqref{eq:HyZ} gives the general form of the asymptotic expansion, including exponentially small contributions.
 
The first exponential correction for any $k\geq 1$ is given by the term $\left(\chi_{11}^{(1)}(x)-\chi_{22}^{(1)}(x)\right)|x|^{1/2}$, because the other two terms contain extra powers of $t^{-1/2}$, either because of the square or from the combination with $Z_{21}^{(1)}(x)$. For $k=1$, we use \eqref{eq:chi111} and 
the leading term is
\[
\begin{aligned}
\left(\left(\chi_{1,1}^{(1)}\right)_{11}
-
\left(\chi_{1,1}^{(1)}\right)_{22}\right)|x|^{1/2}
&=
\frac{2^{-11/8}3^{-1/8}s_{-1}}{\sqrt{\pi}}
\ee^{-\frac{\ii(\varphi-\pi)}{8}}|x|^{-1/8}\ee^{2tg_0}\left(1+\mathcal{O}(x^{-5/4})\right)\\
&=
\frac{2^{-11/8}3^{-1/8}s_{-1}}{\sqrt{\pi}}(\ee^{-\pi\ii}x)^{-1/8}
\ee^{-A(\ee^{-\pi\ii} x)^{5/4}}\left(1+\mathcal{O}(x^{-5/4})\right),
\end{aligned}
\]
with $A=\frac{2^{11/4} 3^{1/4}}{5}$, since from \eqref{eq:glocal_l0} we obtain
\[
2tg_0
=
-\frac{48\sqrt{3}}{5}\lambda_0^{5/2}|x|^{5/4}
=
-\frac{2^{11/4} 3^{1/4}}{5}(\ee^{-\pi\ii}x)^{5/4}.
\]
This gives the first coefficient $y_{1,0}$ in \eqref{eq:y10h10}, and it agrees with \cite[Theorem 2.2]{Kapaev_2004}. Similarly, for the Hamiltonian $\mathcal{H}(x)$ we have the first correction
\[
\begin{aligned}
-\left(\chi_{1,1}^{(1)}(x)\right)_{21}|x|^{1/4}
&=
-\frac{2^{-17/8}3^{-3/8}s_{-1}}{\sqrt{\pi}}
(\ee^{-\pi\ii}x)^{-3/8}\ee^{2tg_0}\left(1+\mathcal{O}(x^{-5/4})\right)\\
&=
-\frac{2^{-17/8}3^{-3/8}s_{-1}}{\sqrt{\pi}}
(\ee^{-\pi\ii}x)^{-3/8}\ee^{-A(\ee^{-\pi\ii} x)^{5/4}}\left(1+\mathcal{O}(x^{-5/4})\right).
\end{aligned}
\]
This gives the first coefficient $h_{1,0}$ in \eqref{eq:y10h10}.

\section*{Acknowledgements}
The author would like to thank Alexander R. Its and Ahmad Barhoumi for very useful discussions on Painlev\'e transcendents and their asymptotics, Silvia N. Santalla for help in elaborating Figure \ref{fig_blob}, and the Isaac Newton Institute for Mathematical Sciences for support and hospitality during the programme Applicable Resurgent Asymptotics: towards a universal theory (ARA), in September--December 2022, when work on this paper was undertaken. This work was supported by EPSRC Grant Number EP/K032208/1. 

The author acknowledges financial support from Universidad Carlos III de Madrid (I Convocatoria para la recualificaci\'on del profesorado universitario), from Direcci\'on General de Investigaci\'on e Innovaci\'on, Consejer\'ia de Educaci\'on e Investigaci\'on of Comunidad de Madrid (Spain), and Universidad de Alcal\'a under grant CM/JIN/2021-014, and Comunidad de Madrid (Spain) under the Multiannual Agreement with UC3M in the line of Excellence of University Professors (EPUC3M23), and in the context of the V PRICIT (Regional Programme of Research and Technological Innovation). Research supported by Grant PID2021-123969NB-I00, funded by MCIN/AEI/ 10.13039/501100011033, and by grant PID2021-122154NB-I00 from Spanish Agencia Estatal de Investigaci\'on.

\bibliographystyle{plain}
\bibliography{biblioPI}

\end{document}